# Bilaplacians problems with a sign-changing coefficient

Lucas Chesnel[1]

[1] Department of Mathematics and Systems Analysis, Aalto University, P.O. Box 11100, FI-00076 Aalto, Finland.

E-mail: `lucas.chesnel@aalto.fi`

**Abstract.** We investigate the properties of the operator $\Delta(\sigma\Delta\cdot) : H^2_0(\Omega) \to H^{-2}(\Omega)$, where $\sigma$ is a given parameter whose sign can change on the bounded domain $\Omega$. Here, $H^2_0(\Omega)$ denotes the subspace of $H^2(\Omega)$ made of the functions $v$ such that $v = \nu \cdot \nabla v = 0$ on $\partial\Omega$. The study of this problem arises when one is interested in some configurations of the Interior Transmission Eigenvalue Problem. We prove that $\Delta(\sigma\Delta\cdot) : H^2_0(\Omega) \to H^{-2}(\Omega)$ is a Fredholm operator of index zero as soon as $\sigma \in L^\infty(\Omega)$, with $\sigma^{-1} \in L^\infty(\Omega)$, is such that $\sigma$ remains uniformly positive (or uniformly negative) in a neighbourhood of $\partial\Omega$. We also study configurations where $\sigma$ changes sign on $\partial\Omega$ and we prove that Fredholm property can be lost for such situations. In the process, we examine in details the features of a simpler problem where the boundary condition $\nu \cdot \nabla v = 0$ is replaced by $\sigma\Delta v = 0$ on $\partial\Omega$.

# Contents

# 1 Introduction

The Interior Transmission Eigenvalue Problem [28, 13] is a spectral problem which appears in inverse scattering theory. In particular, it arises when one is interested in the reconstruction of the support of a penetrable inclusion embedded in a reference medium from far fields measurements at a given frequency. Broadly speaking, the problematic is the following. We consider a reference medium (human body, cable, sea, ...) which contains a localized defect (clot, crack, submarine, ...). Here, "localized" means that the support of the defect is bounded. We send waves in all directions at a given frequency and we measure the scattered field. The goal is to try to deduce from these measurements the support of the defect as well as information on the features of the material of the inclusion (physical parameters).

In this context, it is important to know if for a given frequency, we can find an incident wave for which the field scattered by the inclusion is null. Frequencies for which the answer to this question is positive have an important role: we can use them to characterize the inclusion [8, 22]. There are various objectives in the theory. One consists in proving that the set of these particular frequencies is discrete (see the review paper [14]), another is to establish that such frequencies exist [41, 9, 10]. Here, we propose to consider the first of these two objectives.

In our study, the reference medium will be chosen equal to $\mathbb{R}^d$, $d \geq 1$. The inclusion will be a domain $\Omega \subset \mathbb{R}^d$, *i.e.* a bounded and connected open subset of $\mathbb{R}^d$ with Lipschitz boundary $\partial\Omega$. We assume that the propagation of waves in $\mathbb{R}^d$ is governed by the equation $\Delta w + k^2 w = 0$, where $k \in \mathbb{R}$ is the wave number. On the other hand, we model the inclusion by some physical parameter $n$ so that the total field $u$ (the sum of the incident and scattered fields) satisfies $\Delta u + k^2 n^2 u = 0$. We impose that $n \neq 1$ in $\Omega$ and $n = 1$ in $\mathbb{R}^d \setminus \Omega$. If we denote $[\cdot]|_{\partial\Omega}$ the jump on $\partial\Omega$ (here, the sign is not important) and $\nu$ the unit outward normal vector to $\partial\Omega$ oriented to the exterior of $\Omega$, the total field $u$ satisfies the equations of continuity $[u]|_{\partial\Omega} = [\nu \cdot \nabla u]|_{\partial\Omega} = 0$. Now, if $w$ is an incident field which does not scatter, then there holds $u = w$ outside $\Omega$. As a consequence, we must have $u = w$ and $\nu \cdot \nabla u = \nu \cdot \nabla w$ on $\partial\Omega$. To summarize, if $w$ is such that the scattered field is null outside the inclusion, then the pair $(u, w)$ verifies the problem

$$\left| \begin{array}{lll} \Delta u + k^2 n^2 u & = 0 & \text{in } \Omega \\ \Delta w + k^2 w & = 0 & \text{in } \Omega \\ u - w & = 0 & \text{on } \partial\Omega \\ \nu \cdot \nabla u - \nu \cdot \nabla w & = 0 & \text{on } \partial\Omega. \end{array} \right. \quad (1)$$

Let us introduce some basic notations to equip Problem (1) with a functional framework. The space $\mathrm{L}^2(\Omega)$ is endowed with the classical inner product

$$(\varphi, \varphi')_\Omega = \int_\Omega \varphi \, \overline{\varphi'}, \qquad \forall (\varphi, \varphi') \in \mathrm{L}^2(\Omega) \times \mathrm{L}^2(\Omega).$$

For all $\varphi \in \mathrm{L}^2(\Omega)$, we define $\|\varphi\|_\Omega := (\varphi, \varphi)^{1/2}$. We denote $\mathrm{H}^1_0(\Omega)$ (resp. $\mathrm{H}^2_0(\Omega)$) the closure of $\mathscr{C}^\infty_0(\Omega)$ for the $\mathrm{H}^1$-norm (resp. $\mathrm{H}^2$-norm). We endow these spaces with the inner products

$$\begin{array}{rll} (\varphi, \varphi')_{\mathrm{H}^1_0(\Omega)} & = (\nabla\varphi, \nabla\varphi')_\Omega, & \forall (\varphi, \varphi') \in \mathrm{H}^1_0(\Omega) \times \mathrm{H}^1_0(\Omega); \\ (\varphi, \varphi')_{\mathrm{H}^2_0(\Omega)} & = (\Delta\varphi, \Delta\varphi')_\Omega, & \forall (\varphi, \varphi') \in \mathrm{H}^2_0(\Omega) \times \mathrm{H}^2_0(\Omega). \end{array}$$

The topological dual space of $\mathrm{H}^1_0(\Omega)$ (resp. $\mathrm{H}^2_0(\Omega)$) is denoted $\mathrm{H}^{-1}(\Omega)$ (resp. $\mathrm{H}^{-2}(\Omega)$).

**Definition 1.1.** *The elements $k \in \mathbb{C}$ for which there exists a non trivial solution to the problem*

$$\left| \begin{array}{lll} \textit{Find } (u, w) \in \mathrm{L}^2(\Omega) \times \mathrm{L}^2(\Omega), \textit{ with } u - w \in \mathrm{H}^2_0(\Omega), \textit{ such that:} \\ \Delta u + k^2 n^2 u & = 0 & \textit{in } \Omega \\ \Delta w + k^2 w & = 0 & \textit{in } \Omega \end{array} \right. \quad (2)$$

*are called interior transmission eigenvalues.*

Following [44], to avoid having to work with a system of *PDEs*, we will rewrite (2) as a fourth order equation. Let us describe the method. Consider $(u, w)$ a pair satisfying (2). Define $v := u - w$. Notice that $v$ is nothing else than the scattered field inside the inclusion. It verifies the relation

$$\Delta v + k^2 n^2 v = -k^2(n^2 - 1)w \qquad \text{in } \Omega. \tag{3}$$

In the present paper, we shall assume that the parameter $n : \Omega \to \mathbb{R}$ is an element of $\mathrm{L}^\infty(\Omega)$ such that $n > 0$ and $n \neq 1$ in $\Omega$. We shall also assume that $(n^2 - 1)^{-1}$ belongs to $\mathrm{L}^\infty(\Omega)$[1]. Dividing on each side of (3) by $n^2 - 1$ and using the equation $\Delta w + k^2 w = 0$, we obtain, in the sense of distributions,

$$(\Delta + k^2)\left(\frac{1}{n^2 - 1}(\Delta v + k^2 n^2 v)\right) = 0.$$

We deduce that if the pair $(u, w)$ satisfies Problem (2) then $v = u - w$ verifies the problem

$$\left| \begin{array}{l} \text{Find } v \in \mathrm{H}_0^2(\Omega) \text{ such that} \\ \displaystyle\int_\Omega \frac{1}{n^2 - 1}(\Delta v + k^2 n^2 v)(\overline{\Delta v' + k^2 v'}) = 0, \qquad \forall v' \in \mathrm{H}_0^2(\Omega). \end{array} \right. \tag{4}$$

Conversely, one shows (see [44, lemma 3.1] for the details) that if $v$ is a solution of (4) then the pair $(u, w) := ((n^2 - 1)^{-1}(\Delta v + k^2 n^2 v) - k^2 v, (n^2 - 1)^{-1}(\Delta v + k^2 n^2 v))$ satisfies Problem (2). On $\mathrm{H}_0^2(\Omega) \times \mathrm{H}_0^2(\Omega)$, for $k \in \mathbb{C}$, we define the sesquilinear form $a_k$ such that

$$a_k(v, v') \;=\; ((n^2 - 1)^{-1}(\Delta v + k^2 n^2 v), (\Delta v' + k^2 v'))_\Omega, \qquad \forall (v, v') \in \mathrm{H}_0^2(\Omega) \times \mathrm{H}_0^2(\Omega).$$

With the Riesz representation theorem, let us introduce the bounded operator $A_k : \mathrm{H}_0^2(\Omega) \to \mathrm{H}_0^2(\Omega)$ associated with $a_k$ such that

$$(A_k v, v')_{\mathrm{H}_0^2(\Omega)} = a_k(v, v'), \qquad \forall (v, v') \in \mathrm{H}_0^2(\Omega) \times \mathrm{H}_0^2(\Omega). \tag{5}$$

Assume for a short time that there exists $k_0 \in \mathbb{C}$ such that $A_{k_0}$ is an isomorphism of $\mathrm{H}_0^2(\Omega)$. For all $k \in \mathbb{C}$, it is easy to see that $A_k - A_{k_0}$ is a compact operator of $\mathrm{H}_0^2(\Omega)$. With the analytic Fredholm theorem, we deduce, since $A_{k_0}$ is injective, that $A_k$ is an isomorphism of $\mathrm{H}_0^2(\Omega)$ for all $k \in \mathbb{C} \setminus \mathscr{S}$, where $\mathscr{S}$ is a discrete or empty set of $\mathbb{C}$. In other words, to prove that the set of transmission eigenvalues is discrete, it is sufficient to establish that there exists $k_0 \in \mathbb{C}$ such that $A_{k_0}$ is an isomorphism of $\mathrm{H}_0^2(\Omega)$.

When $n \geq C > 1$ in $\Omega$ for some constant $C$, the bilinear form $a_0$ is coercive on $\mathrm{H}_0^2(\Omega) \times \mathrm{H}_0^2(\Omega)$. This allows to prove that $A_0$ is an isomorphism of $\mathrm{H}_0^2(\Omega)$, and as a consequence, that the set of transmission eigenvalues is discrete or empty. The case $0 < C_1 \leq n \leq C_2 < 1$, where $C_1, C_2$ are some constants, can be dealt with analogously because $-a_0$ is coercive on $\mathrm{H}_0^2(\Omega) \times \mathrm{H}_0^2(\Omega)$. These results are known since [44]. When $n^2 - 1$ changes sign in $\Omega$, the form $a_k$ is no longer coercive nor "coercive+compact" on $\mathrm{H}_0^2(\Omega) \times \mathrm{H}_0^2(\Omega)$. Is the set of transmission eigenvalues still discrete in this case? And before that, what are the properties of $A_k$ in such configurations?

To study the latter question, we will focus on the principal part (the part which contains the derivatives of higher degree) of the operator $A_k$. To simplify the notations, we define $\sigma := (n^2 - 1)^{-1} \in \mathrm{L}^\infty(\Omega)$. All along the paper, we shall assume that $\sigma$ is real valued. Our goal is to investigate the features of the following source term problem

$$(\mathscr{P}) \left| \begin{array}{l} \text{Find } v \in \mathrm{H}_0^2(\Omega) \text{ such that:} \\ (\sigma \Delta v, \Delta v')_\Omega = \langle f, v' \rangle_\Omega, \qquad \forall v' \in \mathrm{H}_0^2(\Omega). \end{array} \right. \tag{6}$$

---

[1] For the study of the Interior Transmission Eigenvalues Problem when $n$ is smooth, we refer the reader to the recent papers [32, 43].

In Problem ($\mathscr{P}$), $f$ is an element of $\mathrm{H}^{-2}(\Omega)$ whereas $\langle \cdot, \cdot \rangle_\Omega$ refers to the duality pairing $\mathrm{H}^{-2}(\Omega) \times \mathrm{H}_0^2(\Omega)$. Let us introduce the sesquilinear form $b$ such that

$$b(v, v') = (\sigma \Delta v, \Delta v')_\Omega, \qquad \forall (v, v') \in \mathrm{H}_0^2(\Omega) \times \mathrm{H}_0^2(\Omega)$$

and the continuous operator $B : \mathrm{H}_0^2(\Omega) \to \mathrm{H}_0^2(\Omega)$ defined by

$$(Bv, v')_{\mathrm{H}_0^2(\Omega)} = b(v, v'), \qquad \forall (v, v') \in \mathrm{H}_0^2(\Omega) \times \mathrm{H}_0^2(\Omega). \tag{7}$$

In recent years, much work (see [15, 6, 2, 4, 40]) has been devoted to the study of the operator $\mathrm{div}(\sigma \nabla \cdot) : \mathrm{H}_0^1(\Omega) \to \mathrm{H}^{-1}(\Omega)$ when $\sigma$ is a parameter whose sign changes on the domain $\Omega$. This operator arises in the modelling of electromagnetic phenomena in time harmonic regime in media involving usual positive material and metals at optical frequencies or negative metamaterials. In this context, the medium is divided into two regions: one corresponding to the positive material ($\sigma = \sigma_+ > 0$), another corresponding to the negative material ($\sigma = \sigma_- < 0$). Let us present the main results, assuming to simplify that $\sigma_+$ and $\sigma_-$ are some constants. If the interface between the two materials is smooth, the operator $\mathrm{div}(\sigma \nabla \cdot) : \mathrm{H}_0^1(\Omega) \to \mathrm{H}^{-1}(\Omega)$ is Fredholm of index zero (see Definition 1.2 below) if and only the contrast $\sigma_+/\sigma_-$ satisfies $\sigma_+/\sigma_- \neq -1$ [15, 2]. When the interface between the two materials has corners, strong singularities can appear. In such configurations, the operator $\mathrm{div}(\sigma \nabla \cdot) : \mathrm{H}_0^1(\Omega) \to \mathrm{H}^{-1}(\Omega)$ is Fredholm of index zero if and only $\sigma_+/\sigma_-$ lies outside some interval $I \subset (-\infty; 0)$ [19, 7, 5, 11]. The latter interval always contains the value $-1$ and depends only on the smallest aperture of the corners of the interface. The goal of the present article is also to compare the features of $\mathrm{div}(\sigma \nabla \cdot) : \mathrm{H}_0^1(\Omega) \to \mathrm{H}^{-1}(\Omega)$ and $\Delta(\sigma \Delta \cdot) : \mathrm{H}_0^2(\Omega) \to \mathrm{H}^{-2}(\Omega)$. Has the change of sign of $\sigma$ the same consequences for both operators?

The outline of the paper is the following. In Section 2, we will study the properties of the operator $\tilde{B} : \mathrm{H}_0^1(\Omega) \cap \mathrm{H}^2(\Omega) \to \mathrm{H}_0^1(\Omega) \cap \mathrm{H}^2(\Omega)$ such that $(\tilde{B}v, v')_{\mathrm{H}_0^2(\Omega)} = (\sigma \Delta v, \Delta v')_\Omega$, for all $v, v'$ in $\mathrm{H}_0^1(\Omega) \cap \mathrm{H}^2(\Omega)$. The latter functional framework corresponds to mixed boundary conditions. We will prove that when the domain $\Omega$ is smooth or convex, $\tilde{B}$ is an isomorphism without assumption on the sign of $\sigma$. Then, adapting a technique used to consider the case $\sigma = 1$ (see [38] and the monograph [20]), we will establish that for polygons with reentrant corners, a kernel and a cokernel can appear for $\tilde{B}$ when $\sigma$ changes sign. The investigation of this simpler problem will provide us the way to study the original operator $B$ (with Dirichlet boundary conditions). This will be the subject of the first part of Section 3 where we will show that $B$ is Fredholm of index zero (see Definition 1.2 hereafter) as soon as $\sigma$ remains uniformly positive or uniformly negative in a neighbourhood of the boundary $\partial\Omega$. In the second part of Section 3, we will be interested in what happens when the sign of $\sigma$ changes on $\partial\Omega$. In particular, we will provide situations where Fredholmness in $\mathrm{H}_0^2(\Omega)$ is lost. Finally in Appendix 4, we will present some features of the operator $\tilde{B}$ in domains with non smooth boundaries in dimension $d \geq 3$, whereas in Appendix 5, we will detail some computations needed in Section 3 in the study of singularities for $B$.

In the sequel, on several occasions, we shall rely on Fredholm theory using the following definition.

**Definition 1.2.** *Let* X *and* Y *be two Banach spaces, and let* $L : \mathrm{X} \to \mathrm{Y}$ *be a continuous linear map. The operator* $L$ *is said to be a Fredholm operator if and only if the following two conditions are fulfilled*

*i)* $\dim(\ker L) < \infty$ *and* range $L$ *is closed;*

*ii)* $\dim(\mathrm{coker}\, L) < \infty$ *where* $\mathrm{coker}\, L := (\mathrm{Y}/\mathrm{range}\, L)$.

*Besides, the index of a Fredholm operator* $L$ *is defined by* $\mathrm{ind}\, L = \dim(\ker L) - \dim(\ker L)$.

## 2 Bilaplacian with mixed boundary conditions

Before investigating the properties of the operator $B$, let us study the problem obtained replacing in $(\mathscr{P})$ the boundary condition "$\nu \cdot \nabla v = 0$ on $\partial\Omega$" by the condition "$\sigma \Delta v = 0$ on $\partial\Omega$".

### 2.1 Smooth or convex domain in $\mathbb{R}^d$, $d \geq 1$

For $f \in (\mathrm{H}_0^1(\Omega) \cap \mathrm{H}^2(\Omega))^*$, the topological dual space of $\mathrm{H}_0^1(\Omega) \cap \mathrm{H}^2(\Omega)$, let us consider the problem

$$\left| \begin{array}{lcl} \text{Find } v \in \mathrm{H}_0^1(\Omega) \cap \mathrm{H}^2(\Omega) & \text{such that:} \\ \Delta(\sigma \Delta v) & = & f \quad \text{in } \Omega \\ \sigma \Delta v & = & 0 \quad \text{on } \partial\Omega. \end{array} \right. \tag{8}$$

Here, we impose mixed boundary conditions: the condition $v = 0$ on $\partial\Omega$ is said to be *essential*, its appears in the functional space, whereas the condition $\sigma \Delta v = 0$ on $\partial\Omega$ is said to be *natural*. The trace $\sigma \Delta v = 0$ is defined in a weak sense. We shall say that the function $\varphi \in \mathrm{L}^2(\Omega)$ such that $\Delta \varphi \in (\mathrm{H}_0^1(\Omega) \cap \mathrm{H}^2(\Omega))^*$ satisfies $\varphi = 0$ on $\partial\Omega$ if and only if there holds

$$\langle \Delta \varphi, \varphi' \rangle_\Omega = (\varphi, \Delta \varphi')_\Omega, \qquad \forall \varphi' \in \mathrm{H}_0^1(\Omega) \cap \mathrm{H}^2(\Omega), \tag{9}$$

where $\langle \cdot, \cdot \rangle_\Omega$ denotes the duality pairing $(\mathrm{H}_0^1(\Omega) \cap \mathrm{H}^2(\Omega))^* \times \mathrm{H}_0^1(\Omega) \cap \mathrm{H}^2(\Omega)$. Therefore, Problem (8) is equivalent to the following problem

$$(\tilde{\mathscr{P}}) \left| \begin{array}{l} \text{Find } v \in \mathrm{H}_0^1(\Omega) \cap \mathrm{H}^2(\Omega) \text{ such that:} \\ (\sigma \Delta v, \Delta v')_\Omega = \langle f, v' \rangle_\Omega, \qquad \forall v' \in \mathrm{H}_0^1(\Omega) \cap \mathrm{H}^2(\Omega). \end{array} \right. \tag{10}$$

With these mixed boundary conditions, it is very easy to solve $(\tilde{\mathscr{P}})$ in two steps. Let $f$ be a source term of $\mathrm{H}^{-1}(\Omega) \subset (\mathrm{H}_0^1(\Omega) \cap \mathrm{H}^2(\Omega))^*$. There exists a unique $p \in \mathrm{H}_0^1(\Omega)$ verifying $-(\nabla p, \nabla p')_\Omega = \langle f, p' \rangle_\Omega$ for all $p' \in \mathrm{H}_0^1(\Omega)$. Let us denote $v$ the unique function satisfying $v \in \mathrm{H}_0^1(\Omega)$ and $\Delta v = \sigma^{-1} p \in \mathrm{L}^2(\Omega)$. If the domain $\Omega$ is of class $\mathscr{C}^2$ ([21, theorem 8.12]) or convex ([23, theorem 3.2.1.2]), we know that $v$ belongs to $\mathrm{H}^2(\Omega)$. Moreover, there holds, for all $v' \in \mathrm{H}_0^1(\Omega) \cap \mathrm{H}^2(\Omega)$,

$$(\sigma \Delta v, \Delta v')_\Omega = (p, \Delta v')_\Omega = -(\nabla p, \nabla v')_\Omega = \langle f, v' \rangle_\Omega.$$

We deduce that $v$ is a solution of $(\tilde{\mathscr{P}})$. Notice that, to obtain this result, the only assumptions for $\sigma$ are $\sigma \in \mathrm{L}^\infty(\Omega)$ and $\sigma^{-1} \in \mathrm{L}^\infty(\Omega)$. Thus, $\sigma$ can change sign. Let us precise the study of $(\tilde{\mathscr{P}})$.

With the Lax-Milgram theorem, we can show that the sesquilinear form $(v, v') \mapsto (v, v)_{\mathrm{H}_0^2(\Omega)} = (\Delta v, \Delta v')_\Omega$ is an inner product on $\mathrm{H}_0^1(\Omega) \cap \mathrm{H}^2(\Omega)$. Moreover, if $\Omega$ is of class $\mathscr{C}^2$ or convex, on $\mathrm{H}_0^1(\Omega) \cap \mathrm{H}^2(\Omega)$, the map $v \mapsto \|\Delta v\|_\Omega$ defines a norm which is equivalent to the $\mathrm{H}^2$-norm. Therefore, $\mathrm{H}_0^1(\Omega) \cap \mathrm{H}^2(\Omega)$ endowed with the inner product $(\cdot, \cdot)_{\mathrm{H}_0^2(\Omega)}$ is a Hilbert space. Let us introduce the sesquilinear form $\tilde{b}$ such that

$$\tilde{b}(v, v') = (\sigma \Delta v, \Delta v')_\Omega, \qquad \forall (v, v') \in \mathrm{H}_0^1(\Omega) \cap \mathrm{H}^2(\Omega) \times \mathrm{H}_0^1(\Omega) \cap \mathrm{H}^2(\Omega), \tag{11}$$

and the continuous operator $\tilde{B} : \mathrm{H}_0^1(\Omega) \cap \mathrm{H}^2(\Omega) \to \mathrm{H}_0^1(\Omega) \cap \mathrm{H}^2(\Omega)$, defined by

$$(\tilde{B} v, v')_{\mathrm{H}_0^2(\Omega)} = \tilde{b}(v, v'), \qquad \forall (v, v') \in \mathrm{H}_0^1(\Omega) \cap \mathrm{H}^2(\Omega) \times \mathrm{H}_0^1(\Omega) \cap \mathrm{H}^2(\Omega). \tag{12}$$

**Theorem 2.1.** *Assume that the domain $\Omega \subset \mathbb{R}^d$, $d \geq 1$, is of class $\mathscr{C}^2$ or convex. For all $\sigma \in \mathrm{L}^\infty(\Omega)$ such that $\sigma^{-1} \in \mathrm{L}^\infty(\Omega)$, the operator $\tilde{B} : \mathrm{H}_0^1(\Omega) \cap \mathrm{H}^2(\Omega) \to \mathrm{H}_0^1(\Omega) \cap \mathrm{H}^2(\Omega)$ defined in (12) is an isomorphism.*

*Proof.* Let us introduce the operator $\mathtt{T} : \mathrm{H}^1_0(\Omega) \cap \mathrm{H}^2(\Omega) \to \mathrm{H}^1_0(\Omega) \cap \mathrm{H}^2(\Omega)$ such that, for all $v \in \mathrm{H}^1_0(\Omega) \cap \mathrm{H}^2(\Omega)$, $\mathtt{T}v$ is defined as the unique solution of the problem "find $\mathtt{T}v \in \mathrm{H}^1_0(\Omega)$ satisfying $\Delta(\mathtt{T}v) = \sigma^{-1}\Delta v$". Notice that since $\Omega$ is assumed to be of class $\mathscr{C}^2$ ([21, theorem 8.12]) or convex ([23, theorem 3.2.1.2]), $\mathtt{T}v$ is indeed an element of $\mathrm{H}^2(\Omega)$. For all $v, v' \in \mathrm{H}^1_0(\Omega) \cap \mathrm{H}^2(\Omega)$, we can write
$$(\tilde{B}(\mathtt{T}v), v')_{\mathrm{H}^2_0(\Omega)} = (\sigma \Delta(\mathtt{T}v), \Delta v')_\Omega = (\Delta v, \Delta v')_\Omega.$$

Therefore, the operator $\tilde{B} \circ \mathtt{T}$ is equal to the identity of $\mathrm{H}^1_0(\Omega) \cap \mathrm{H}^2(\Omega)$. Since $\tilde{B}$ is selfadjoint, we deduce that $\tilde{B}$ is an isomorphism with $\tilde{B}^{-1} = \mathtt{T}$. □

Let us give another proof of this result, slightly different, using the resolution of $(\tilde{\mathscr{P}})$ in two steps. This will give us an idea of how to proceed to study configurations where $\partial\Omega$ is not smooth. Let us start by proving the classical result

**Proposition 2.1.** *Assume that the domain $\Omega \subset \mathbb{R}^d$, $d \geq 1$, is of class $\mathscr{C}^2$ or convex. Then for all $f \in (\mathrm{H}^1_0(\Omega) \cap \mathrm{H}^2(\Omega))^*$, there exists a unique solution to the problem*

$$\left| \begin{array}{l} \textit{Find } p \in \mathrm{L}^2(\Omega) \textit{ such that:} \\ (p, \Delta v')_\Omega = \langle f, v' \rangle_\Omega, \qquad \forall v' \in \mathrm{H}^1_0(\Omega) \cap \mathrm{H}^2(\Omega) \end{array} \right. \tag{13}$$

*Here, $\langle \cdot, \cdot \rangle_\Omega$ refers to the duality pairing $(\mathrm{H}^1_0(\Omega) \cap \mathrm{H}^2(\Omega))^* \times \mathrm{H}^1_0(\Omega) \cap \mathrm{H}^2(\Omega)$.*

*Proof.* Since $\Omega$ is smooth or convex, $\Delta : \mathrm{H}^1_0(\Omega) \cap \mathrm{H}^2(\Omega) \to \mathrm{L}^2(\Omega)$ is an isomorphism. The result of this proposition is obtained thanks to the well-known proposition which asserts that the adjoint operator of an isomorphism is also an isomorphism. Let us detail the proof of this point in this particular case.

Let $p$ be an element of $\mathrm{L}^2(\Omega)$ such that $(p, \Delta v')_\Omega = 0$ for all $v' \in \mathrm{H}^1_0(\Omega) \cap \mathrm{H}^2(\Omega)$. Testing with $v'$ the unique solution to the problem "find $v' \in \mathrm{H}^1_0(\Omega) \cap \mathrm{H}^2(\Omega)$ such that $\Delta v' = p$", we deduce that $\|p\|^2_\Omega = 0$. Thus Problem (13) admits at most one solution. Now, consider $f \in (\mathrm{H}^1_0(\Omega) \cap \mathrm{H}^2(\Omega))^*$. According to the Riesz representation theorem, there exists a unique $F \in \mathrm{H}^1_0(\Omega) \cap \mathrm{H}^2(\Omega)$ such that $(\Delta F, \Delta v')_\Omega = \langle f, v' \rangle_\Omega$ for all $v' \in \mathrm{H}^1_0(\Omega) \cap \mathrm{H}^2(\Omega)$. The function $p = \Delta F \in \mathrm{L}^2(\Omega)$ is then a solution to Problem (13). This ends the proof. □

*Proof of Theorem 2.1 (bis).* Let $v \in \mathrm{H}^1_0(\Omega) \cap \mathrm{H}^2(\Omega)$ such that
$$(\sigma \Delta v, \Delta v')_\Omega = 0, \qquad \forall v' \in \mathrm{H}^1_0(\Omega) \cap \mathrm{H}^2(\Omega).$$

According to Proposition 2.1, we have necessarily $\sigma \Delta v = 0$, hence $\Delta v = 0$. This implies $v = 0$ since $v \in \mathrm{H}^1_0(\Omega) \cap \mathrm{H}^2(\Omega)$ and proves that $(\tilde{\mathscr{P}})$ admits at most one solution. Then, we consider a source term $f$ in $(\mathrm{H}^1_0(\Omega) \cap \mathrm{H}^2(\Omega))^*$. By virtue of Proposition 2.1, there exists a unique $p \in \mathrm{L}^2(\Omega)$ such that
$$(p, \Delta v')_\Omega = \langle f, v' \rangle_\Omega, \qquad \forall v' \in \mathrm{H}^1_0(\Omega) \cap \mathrm{H}^2(\Omega).$$

Let us call $v$ the unique solution to the problem "find $v \in \mathrm{H}^1_0(\Omega) \cap \mathrm{H}^2(\Omega)$ such that $\Delta v = \sigma^{-1}p$". This function $v$ is a solution to $(\tilde{\mathscr{P}})$. □

**Remark 2.1.** *The resolution of $(\tilde{\mathscr{P}})$ in two steps is interesting for numerical considerations. Indeed, it is easier to solve two second order problems than a fourth order one.*

The operator $\tilde{B}$ defined in (12) is an isomorphism of $\mathrm{H}^1_0(\Omega) \cap \mathrm{H}^2(\Omega)$ when $\sigma \in \mathrm{L}^\infty(\Omega)$ verifies $\sigma^{-1} \in \mathrm{L}^\infty(\Omega)$ and when $\Omega$ is such that $\Delta : \mathrm{H}^1_0(\Omega) \cap \mathrm{H}^2(\Omega) \to \mathrm{L}^2(\Omega)$ constitutes an isomorphism. But what if this last assumption on $\Omega$ is not met? What happens for example if $\Omega$ is a 2D domain with reentrant corners? We prove in the next section that in this case, $\tilde{B}$ is not always an isomorphism. According to the values of the parameter $\sigma$, a kernel and a cokernel, both of finite dimension, can appear. For results in $\mathbb{R}^d$ with $d \geq 3$, we refer the reader to Section 4.

## 2.2 Polygonal domains with reentrant corners

We assume in this paragraph that $\Omega \subset \mathbb{R}^2$ is an open set with a polygonal boundary $\partial\Omega$. For such domains, using an integration by parts, one proves the *a priori* estimate (see [24, theorem 2.2.3] or [26, 34]):

$$\|v\|_{H^2(\Omega)} \leq C \, \|\Delta v\|_{\Omega}, \qquad \forall v \in H_0^1(\Omega) \cap H^2(\Omega),$$

where $C$ is a constant which depends only on $\Omega$. This estimate provides lot of information. It proves that the operator $\Delta : H_0^1(\Omega) \cap H^2(\Omega) \to L^2(\Omega)$ is injective and that its range is closed (it is a monomorphism). Then we can try to characterize the orthogonal complement of the range of $\Delta : H_0^1(\Omega) \cap H^2(\Omega) \to L^2(\Omega)$. Theorem 2.3.7 of [24] indicates that this orthogonal complement is of finite dimension $N$, where $N$ is equal to the number of corners of $\partial\Omega$ whose aperture is strictly larger than $\pi$. Thus, $\Delta : H_0^1(\Omega) \cap H^2(\Omega) \to L^2(\Omega)$ is an injective Fredholm operator of index $-N$. When there is no reentrant corner in $\partial\Omega$, *i.e.* when $\Omega$ is convex, we find back that $\Delta : H_0^1(\Omega) \cap H^2(\Omega) \to L^2(\Omega)$ is an isomorphism.

### 2.2.1 Sapongyan paradox in the case of a positive $\sigma$

We assume in this paragraph that there exists a constant $C$ such that $\sigma \geq C > 0$ a.e. in $\Omega$. In this case, the form $\tilde{b}$ defined in (11) is coercive on $H_0^1(\Omega) \cap H^2(\Omega) \times H_0^1(\Omega) \cap H^2(\Omega)$ whether the domain $\Omega$ is convex or not. According to the Lax-Milgram theorem, $(\tilde{\mathscr{P}})$ has a unique solution $v \in H_0^1(\Omega) \cap H^2(\Omega)$ and $\tilde{B}$ is an isomorphism of $H_0^1(\Omega) \cap H^2(\Omega)$.

Now, let us try to solve $(\tilde{\mathscr{P}})$ in two steps. To simplify the explanation, we assume here that $f$ belongs to $H^{-1}(\Omega)$. Let us denote $p_0 \in H_0^1(\Omega)$ the function such that $-(\nabla p_0, \nabla p')_\Omega = \langle f, p' \rangle_\Omega$ for all $p' \in H_0^1(\Omega)$. Then, let us introduce $v_0$ the element of $H_0^1(\Omega)$ satisfying $\Delta v_0 = \sigma^{-1} p_0 \in L^2(\Omega)$. When $\Omega$ is convex, the function $v_0$ is an element of $H_0^1(\Omega) \cap H^2(\Omega)$. In this case, $v_0$ verifies $(\tilde{\mathscr{P}})$. Since this problem is well-posed, we deduce $v_0 = v$. When $\partial\Omega$ has one or several reentrant corners, it can happen that $v_0 \notin H_0^1(\Omega) \cap H^2(\Omega)$. In this situation, we have $(\sigma \Delta v_0, \Delta v')_\Omega = \langle f, v' \rangle_\Omega$ for all $v' \in H_0^1(\Omega) \cap H^2(\Omega)$ but $v_0$ is not equal to the solution $v$ of $(\tilde{\mathscr{P}})$. This is what S.A. Nazarov and G.H. Sweers call, in the very interesting papers [37, 38, 39], the *Sapongyan's paradox*. Sapongyan was a Russian mathematician from the early XIXth century. Thanks to conformal mapping techniques, he obtained a solution which was not of finite mechanical energy. Since he had no explanation for that, he called this phenomenon a paradox.

Nevertheless, and this is the topic of papers [37, 38, 39], even when $\Omega$ has a reentrant corner, there exists a mean of solving $(\tilde{\mathscr{P}})$ in two steps while obtaining the solution of finite energy, that is the one in $H_0^1(\Omega) \cap H^2(\Omega)$. Let us describe the process. To simplify, we assume that $\Omega$ has only one reentrant corner $O$ of aperture equal to $\alpha \in (\pi; 2\pi)$ (cf. Figure 1). The method consists in solving in a clever way the two Laplace problems with homogeneous Dirichlet boundary condition which appears in $(\tilde{\mathscr{P}})$. Before proceeding further, we need to recall some classical results (see for example [23, 24]) of the theory of singularities for the Laplace operator in non convex polygons.

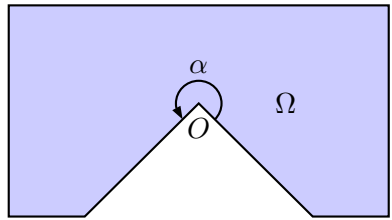

Figure 1: An example of a non convex polygonal boundary with one reentrant corner – $N = 1$.

Let us introduce $\zeta$ such that

$$\zeta(\boldsymbol{x}) = r^{-\pi/\alpha} \sin(\pi\theta/\alpha) + \tilde{\zeta}(\boldsymbol{x}). \tag{14}$$

where $\tilde{\zeta}$ is the unique function of $\mathrm{H}^1(\Omega)$ satisfying $\Delta\tilde{\zeta} = 0$ a.e. in $\Omega$ and $\tilde{\zeta} = -r^{-\pi/\alpha}\sin(\pi\theta/\alpha)$ a.e. on $\partial\Omega$. Here, $(r,\theta)$ are the polar coordinates centered at $O$, such that $\theta = 0$ or $\theta = \alpha$ on $\partial\Omega$ in a neighbourhood of $O$. We assume that $\Omega$ is not convex at $O$. This imposes $0 < \pi/\alpha < 1$. By a straightforward computation, one proves that $\boldsymbol{x} \mapsto r^{-\pi/\alpha}\sin(\pi\theta/\alpha)$ belongs to $\mathrm{L}^2(\Omega)\setminus\mathrm{H}^1(\Omega)$. Since $\tilde{\zeta} \in \mathrm{H}^1(\Omega)$, we deduce that the function $\zeta$ defined in (14) satisfies $\zeta \in \mathrm{L}^2(\Omega)\setminus\mathrm{H}^1(\Omega)$. By definition of $\tilde{\zeta}$, there hold $\Delta\zeta = 0$ a.e. in $\Omega$ and $\zeta = 0$ a.e. on $\partial\Omega$. According to (9), this is equivalent to the following integral identity: $(\zeta, \Delta v')_\Omega = 0$ for all $v' \in \mathrm{H}_0^1(\Omega) \cap \mathrm{H}^2(\Omega)$. Actually, $\zeta$ constitutes a basis of the set of $\mathrm{L}^2$-functions satisfying such property (see [24, lemma 2.3.6]). Let us include this result in the following proposition.

**Proposition 2.2.** *Assume that the boundary $\partial\Omega$ is a polygon which has exactly one reentrant corner. Then for all $f \in (\mathrm{H}_0^1(\Omega) \cap \mathrm{H}^2(\Omega))^*$, there exists a solution to the problem*

$$\left|\begin{array}{l} \textit{Find } p \in \mathrm{L}^2(\Omega) \textit{ such that:} \\ (p, \Delta v')_\Omega = \langle f, v'\rangle_\Omega, \qquad \forall v' \in \mathrm{H}_0^1(\Omega) \cap \mathrm{H}^2(\Omega). \end{array}\right. \tag{15}$$

*Moreover, if $p_1$, $p_2$ are two solutions of Problem (15), then there exists $a \in \mathbb{C}$ such that $p_2 = p_1 + a\zeta$, where $\zeta$ is introduced in (14).*

*Proof.* One establishes that, for all $f \in (\mathrm{H}_0^1(\Omega) \cap \mathrm{H}^2(\Omega))^*$, Problem (15) has at least one solution working as in the proof of Proposition 2.1. Now, if $p_1$, $p_2$ are two solutions of Problem (15), then we have $(p_2 - p_1, \Delta v')_\Omega = 0$ for all $v' \in \mathrm{H}_0^1(\Omega) \cap \mathrm{H}^2(\Omega)$. With [24, lemma 2.3.6], we deduce that there holds $p_2 - p_1 = a\zeta$, where $a$ is a constant. □

Now, we provide a result of decomposition (see for example [24, theorem 2.4.3] or [30, p. 263]) which states that the elements of $\mathrm{H}_0^1(\Omega)$ whose Laplacian belongs to $\mathrm{L}^2(\Omega)$ split as the sum of an explicit singular part and a regular part in $\mathrm{H}^2(\Omega)$.

**Proposition 2.3.** *Consider $\varphi \in \mathrm{H}_0^1(\Omega)$ such that $\Delta\varphi = g \in \mathrm{L}^2(\Omega)$. Then $\varphi$ admits the decomposition*

$$\varphi(\boldsymbol{x}) = c\, r^{\pi/\alpha}\sin(\pi\theta/\alpha) + \tilde{\varphi}(\boldsymbol{x}), \tag{16}$$

*with $\tilde{\varphi} \in \mathrm{H}^2(\Omega)$. Moreover, the coefficient $c$ in (16) is given by the following expression*

$$c = -(\pi)^{-1}(g, \zeta)_\Omega,$$

*where $\zeta$ is defined in (14).*

Noticing that $\boldsymbol{x} \mapsto r^{\pi/\alpha}\sin(\pi\theta/\alpha)$ belongs to $\mathrm{H}^1(\Omega)\setminus\mathrm{H}^2(\Omega)$ (since $\pi/\alpha < 1$), we deduce the

**Corollary 2.1.** *Let $\varphi \in \mathrm{H}_0^1(\Omega)$ be such that $\Delta\varphi = g \in \mathrm{L}^2(\Omega)$. Then $\varphi \in \mathrm{H}^2(\Omega)$ if and only if $(g, \zeta)_\Omega = 0$.*

Now, we have all the tools we need to solve $(\tilde{\mathscr{P}})$ in two steps. For any source term $f \in (\mathrm{H}_0^1(\Omega) \cap \mathrm{H}^2(\Omega))^*$, let us introduce $p_0$ an element of $\mathrm{L}^2(\Omega)$ such that $(p_0, \Delta v')_\Omega = \langle f, v'\rangle_\Omega$ for all $v' \in \mathrm{H}_0^1(\Omega) \cap \mathrm{H}^2(\Omega)$. The existence of such a function $p_0$ is guaranteed by Proposition 2.2. Now, for $a \in \mathbb{C}$, we define $p = p_0 + a\zeta$, where $\zeta$ is given by (14), and denote $v$ the unique function of $\mathrm{H}_0^1(\Omega)$ such that $\Delta v = \sigma^{-1}p$. We want $v$ to be in $\mathrm{H}^2(\Omega)$. According to Corollary 2.1, we must take $a$ such that

$$\begin{aligned} 0 &= (\sigma^{-1}p, \zeta)_\Omega = (\sigma^{-1}p_0, \zeta)_\Omega + a\,(\sigma^{-1}\zeta, \zeta)_\Omega \\ \Leftrightarrow \quad a &= -(\sigma^{-1}p_0, \zeta)_\Omega/(\sigma^{-1}\zeta, \zeta)_\Omega. \end{aligned} \tag{17}$$

To conclude that $v$ constitutes the solution of $(\tilde{\mathscr{P}})$, it just remains to notice that, for all $v' \in \mathrm{H}_0^1(\Omega) \cap \mathrm{H}^2(\Omega)$, there holds

$$(\sigma\Delta v, \Delta v')_\Omega = (p, \Delta v')_\Omega = (p_0, \Delta v')_\Omega = \langle f, v'\rangle_\Omega.$$

### 2.2.2 Study in the case where σ changes sign

When $\sigma$ changes sign, the sesquilinear form $\tilde{b}$ associated with $(\tilde{\mathscr{P}})$ is not coercive. When the domain $\Omega$ is convex or of class $\mathscr{C}^2$, constructing the inverse of $\tilde{B}$, we proved with Theorem 2.1 that $(\tilde{\mathscr{P}})$ is well-posed. Moreover, in the second proof of Theorem 2.1, we established that under one of these two assumptions, $(\tilde{\mathscr{P}})$ can also be solved in two steps. In the previous paragraph, for a positive $\sigma$, we presented how to solve $(\tilde{\mathscr{P}})$ in two steps when $\Omega$ has one reentrant corner. Our goal is to extend this approach to deal with configurations where $\sigma$ changes sign. The novelty is that, according to the values of $\sigma$, a kernel and a cokernel, whose dimensions are less or equal to the number of reentrant corners of the domain, can appear.

*2.2.2.1 Polygonal domains with one reentrant corner*

To begin with, we assume that $\partial\Omega$ has only one reentrant corner located at $O$, of aperture $\alpha \in (\pi; 2\pi)$. Working as in §2.2.1, we see that the resolution in two steps can be used to prove that the operator $\tilde{B}$ associated with $(\tilde{\mathscr{P}})$ is onto as soon as $\sigma$ satisfies $(\sigma^{-1}\zeta, \zeta)_\Omega \neq 0$. This leads us to consider two cases: either $\sigma$ is such that $(\sigma^{-1}\zeta, \zeta)_\Omega \neq 0$ or $\sigma$ is such that $(\sigma^{-1}\zeta, \zeta)_\Omega = 0$.

$$\star\star\star\star\star$$
$$\text{Case } (\sigma^{-1}\zeta, \zeta)_\Omega \neq 0$$
$$\star\star\star\star\star$$

Let us assume first that $(\sigma^{-1}\zeta, \zeta)_\Omega \neq 0$. We have the

**Proposition 2.4.** *Assume that the boundary $\partial\Omega$ is a polygon which has exactly one reentrant corner. Assume that $\sigma \in L^\infty(\Omega)$ verifies $\sigma^{-1} \in L^\infty(\Omega)$ and $(\sigma^{-1}\zeta, \zeta)_\Omega \neq 0$, where $\zeta$ is defined in (14). Then, the operator $\tilde{B}: H_0^1(\Omega) \cap H^2(\Omega) \to H_0^1(\Omega) \cap H^2(\Omega)$ defined in (12) is an isomorphism.*

*Proof.* Let us introduce the operator $\mathtt{T}$ such that, for $v \in H_0^1(\Omega) \cap H^2(\Omega)$, the function $\mathtt{T}v \in H_0^1(\Omega)$ satisfies $\Delta(\mathtt{T}v) = \sigma^{-1}(\Delta v + a\zeta)$ with $a = -(\sigma^{-1}\Delta v, \zeta)_\Omega / (\sigma^{-1}\zeta, \zeta)_\Omega$. Since $(\sigma^{-1}(\Delta v + a\zeta), \zeta)_\Omega = 0$, we know that $\mathtt{T}v$ belongs to $H^2(\Omega)$ according to Corollary 2.1. Thus, $\mathtt{T}$ is a continuous operator from $H_0^1(\Omega) \cap H^2(\Omega)$ to $H_0^1(\Omega) \cap H^2(\Omega)$. For all $(v, v') \in H_0^1(\Omega) \cap H^2(\Omega) \times H_0^1(\Omega) \cap H^2(\Omega)$, we then compute

$$\begin{aligned}(\tilde{B}(\mathtt{T}v), v')_{H_0^2(\Omega)} &= \tilde{b}(\mathtt{T}v, v') &= (\sigma\Delta(\mathtt{T}v), \Delta v')_\Omega \\ & &= (\Delta v + a\zeta, \Delta v')_\Omega \\ & &= (\Delta v, \Delta v')_\Omega = (v, v')_{H_0^2(\Omega)}.\end{aligned}$$

The last line is obtained noticing that $(\zeta, \Delta v')_\Omega = 0$ because $v' \in H_0^1(\Omega) \cap H^2(\Omega)$ (Corollary 2.1). Thus, there holds $\tilde{B} \circ \mathtt{T} = \mathrm{Id}$. Since $\tilde{B}$ is selfadjoint, we deduce that $\tilde{B}$ is an isomorphism with $\tilde{B}^{-1} = \mathtt{T}$. □

**Remark 2.2.** *The method developed to prove well-posedness of Problem $(\tilde{\mathscr{P}})$ in geometries with reentrant corners presents strong analogies with the one used to establish well-posedness of Maxwell's equations when the physical parameters $\varepsilon$, $\mu$ change sign on the domain. For the investigation of the latter problem, we refer the reader to [3].*

Let us stop here the study of the case $(\sigma^{-1}\zeta, \zeta)_\Omega \neq 0$ and let us focus our attention on the configuration $(\sigma^{-1}\zeta, \zeta)_\Omega = 0$.

$$\star\star\star\star\star$$
$$\text{Case } (\sigma^{-1}\zeta, \zeta)_\Omega = 0$$
$$\star\star\star\star\star$$

In this case, we can no longer use the precious degree of freedom to construct a solution to $(\tilde{\mathscr{P}})$ in $H_0^1(\Omega) \cap H^2(\Omega)$. Let us denote $\psi$ the function of $H_0^1(\Omega)$ satisfying

$$\Delta\psi = \sigma^{-1}\zeta. \tag{18}$$

Since $(\sigma^{-1}\zeta,\zeta)_\Omega = 0$, Corollary 2.1 indicates that $\psi$ belongs to $\mathrm{H}^2(\Omega)$. Moreover, for all $v' \in \mathrm{H}^1_0(\Omega) \cap \mathrm{H}^2(\Omega)$, there holds $(\sigma\Delta\psi, \Delta v')_\Omega = (\zeta, \Delta v')_\Omega = 0$. Therefore, $\psi$ constitutes an element of $\ker \tilde{B}$.

**Proposition 2.5.** *Assume that the boundary $\partial\Omega$ is a polygon which has exactly one reentrant corner. Assume that $\sigma \in \mathrm{L}^\infty(\Omega)$ verifies $\sigma^{-1} \in \mathrm{L}^\infty(\Omega)$ and $(\sigma^{-1}\zeta,\zeta)_\Omega = 0$. Then,*
- $\dim(\ker \tilde{B}) = 1$ *with* $\ker \tilde{B} = \mathrm{span}(\psi)$;
- $\dim(\mathrm{coker}\,\tilde{B}) = 1$ *and for all $f \in (\mathrm{H}^1_0(\Omega) \cap \mathrm{H}^2(\Omega))^*$, $(\tilde{\mathscr{P}})$ has a solution if and only if $\langle f, \psi \rangle_\Omega = 0$.*

*In this statement, the functions $\zeta$ and $\psi$ are respectively defined in (14) and (18).*

*Proof.* ⋆ Kernel. If $v$ belongs to $\ker \tilde{B}$ then, according to Proposition 2.2, we have $\sigma\Delta v = a\zeta$ where $a$ is a constant. Thus, $\ker \tilde{B} \subset \mathrm{span}(\psi)$. As indicated above, we have $\psi \in \ker \tilde{B}$ and so $\ker \tilde{B} = \mathrm{span}(\psi)$.

⋆ Cokernel. Let us consider $f \in (\mathrm{H}^1_0(\Omega) \cap \mathrm{H}^2(\Omega))^*$ such that $\langle f, \psi \rangle_\Omega = 0$. Let us introduce $p_0$ an element of $\mathrm{L}^2(\Omega)$ such that $(p_0, \Delta v')_\Omega = \langle f, v' \rangle_\Omega$ for all $v' \in \mathrm{H}^1_0(\Omega) \cap \mathrm{H}^2(\Omega)$. The existence of such a function $p_0$ is ensured by Proposition 2.2. The function $p_0$ satisfies the compatibility condition $(\sigma^{-1} p_0, \zeta)_\Omega = 0$. Indeed, since $\Delta\psi = \sigma^{-1}\zeta$, we can write $(\sigma^{-1} p_0, \zeta)_\Omega = (p_0, \sigma^{-1}\zeta)_\Omega = (p_0, \Delta\psi)_\Omega = \langle f, \psi \rangle_\Omega = 0$. As a consequence, by virtue of Corollary 2.1, the function $v \in \mathrm{H}^1_0(\Omega)$ verifying $\Delta v = \sigma^{-1} p_0$ is in $\mathrm{H}^2(\Omega)$. Moreover, for all $v' \in \mathrm{H}^1_0(\Omega) \cap \mathrm{H}^2(\Omega)$, there holds

$$(\sigma\Delta v, \Delta v')_\Omega = (p_0, \Delta v')_\Omega = \langle f, v' \rangle_\Omega.$$

Therefore, $v$ constitutes a solution to $(\tilde{\mathscr{P}})$. Now, let us consider $f \in (\mathrm{H}^1_0(\Omega) \cap \mathrm{H}^2(\Omega))^*$ such that $\langle f, \psi \rangle_\Omega \neq 0$. Let us assume that there exists a solution $v$ to $(\tilde{\mathscr{P}})$. Then, by Proposition 2.2, we have $\sigma\Delta v = p_0 + a\zeta$, where $p_0 \in \mathrm{L}^2(\Omega)$ satisfies $(p_0, \Delta v')_\Omega = \langle f, v' \rangle_\Omega$ for all $v' \in \mathrm{H}^1_0(\Omega) \cap \mathrm{H}^2(\Omega)$ and where $a$ is a constant. This imposes $(\sigma^{-1} p_0, \zeta)_\Omega = 0$. But there holds $(\sigma^{-1} p_0, \zeta)_\Omega = (p_0, \sigma^{-1}\zeta)_\Omega = (p_0, \Delta\psi)_\Omega = \langle f, \psi \rangle_\Omega$. Thus, we obtain an absurdity. This ends to prove that there exists a solution to $(\tilde{\mathscr{P}})$ if and only if $\langle f, \psi \rangle_\Omega = 0$. □

*2.2.2.2 Polygonal domains with several reentrant corners*

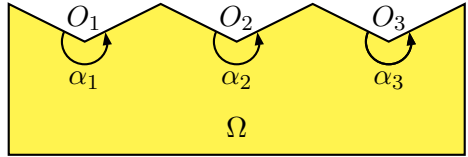

Figure 2: An example of a non convex polygonal boundary with three reentrant corners – $N = 3$.

Assume that the boundary $\partial\Omega$ has $N$ reentrant corners $O_i$ of aperture $\alpha_i \in (\pi; 2\pi)$, $i = 1 \ldots N$ (see Figure 2). We denote $(r_i, \theta_i)$ the polar coordinates centered at $O_i$, such that $\theta_i = 0$ or $\theta_i = \alpha_i$ on $\partial\Omega$ in a neighbourhood of $O_i$. Let us start by recalling some classical results concerning the theory of singularities for the Laplace operator with Dirichlet boundary condition in polygonal domains with several reentrant corners.

For $i = 1 \ldots N$, we introduce $\zeta_i \in \mathrm{L}^2(\Omega) \setminus \mathrm{H}^1(\Omega)$ the function such that

$$\zeta_i(\boldsymbol{x}) = r_i^{-\pi/\alpha_i} \sin(\pi\theta_i/\alpha_i) + \tilde{\zeta}_i(\boldsymbol{x}), \tag{19}$$

where $\tilde{\zeta}_i$ is the unique element of $\mathrm{H}^1(\Omega)$ verifying $\Delta\tilde{\zeta}_i = 0$ a.e. in $\Omega$ and $\tilde{\zeta}_i = -r_i^{-\pi/\alpha_i} \sin(\pi\theta_i/\alpha_i)$ a.e. on $\partial\Omega$. Notice that $\Delta\zeta_i = 0$ a.e. in $\Omega$ and $\zeta_i = 0$ a.e. on $\partial\Omega$. The following result extends the one of Proposition 2.2 to the case where the boundary $\partial\Omega$ contains several reentrant corners.

**Proposition 2.6.** *Assume that the boundary $\partial\Omega$ is a polygon which has exactly $N$ reentrant corners. Then for all $f \in (\mathrm{H}_0^1(\Omega) \cap \mathrm{H}^2(\Omega))^*$, there exists a solution to the problem*

$$\left|\begin{array}{l} \text{Find } p \in \mathrm{L}^2(\Omega) \text{ such that:} \\ (p, \Delta v')_\Omega = \langle f, v' \rangle_\Omega, \qquad \forall v' \in \mathrm{H}_0^1(\Omega) \cap \mathrm{H}^2(\Omega). \end{array}\right. \tag{20}$$

*Moreover, the family $(\zeta_i)_{i=1}^N$, where the $\zeta_i$ are defined in (19), constitutes a basis of the set of functions satisfying Problem (20) with $f = 0$. Therefore, if $p_1$, $p_2$ are two solutions of Problem (20), then there exists some constants $a_1, \ldots, a_N$ such that $p_2 = p_1 + \sum_{i=1}^N a_i \zeta_i$.*

*Proof.* Proceeding as in the proof of Proposition 2.1, one establishes that, for all $f \in (\mathrm{H}_0^1(\Omega) \cap \mathrm{H}^2(\Omega))^*$, Problem (20) has at least one solution. Now, consider $p_1$, $p_2$ two solutions of Problem (20). We have $(p_2 - p_1, \Delta v')_\Omega = 0$ for all $v' \in \mathrm{H}_0^1(\Omega) \cap \mathrm{H}^2(\Omega)$. Thanks to [24, lemma 2.3.6], we deduce that $p_2 = p_1 + \sum_{i=1}^N a_i \zeta_i$ where $a_1, \ldots, a_N$ are some constants. On the other hand, if $a_1, \ldots, a_N$ are some constants such that $\sum_{i=1}^N a_i \zeta_i = 0$ a.e. in $\Omega$, since the functions $\boldsymbol{x} \mapsto r_i^{-\pi/\alpha_i} \sin(\pi\theta_i/\alpha_i)$ are not in $\mathrm{H}^1$ only locally in a neighbourhood of $O_i$, we deduce that $a_1 = \cdots = a_N = 0$. This ends to prove that $(\zeta_i)_{i=1}^N$ constitutes a basis of the set of functions satisfying Problem (20) with $f = 0$. □

In presence of $N$ reentrant corners in the boundary, the result of decomposition which states that the elements of $\mathrm{H}_0^1(\Omega)$ whose Laplacian belongs to $\mathrm{L}^2(\Omega)$ split as the sum of an explicit singular part and a regular part in $\mathrm{H}^2(\Omega)$, becomes (cf. theorem 6.4.4 and paragraph 6.6.1. of [30]):

**Proposition 2.7.** *Consider $\varphi \in \mathrm{H}_0^1(\Omega)$ such that $\Delta\varphi = g \in \mathrm{L}^2(\Omega)$. Then $\varphi$ admits the decomposition*

$$\varphi(\boldsymbol{x}) = \sum_{i=1}^N c_i \, r_i^{\pi/\alpha_i} \sin(\pi\theta_i/\alpha_i) + \tilde{\varphi}_i(\boldsymbol{x}), \tag{21}$$

*with $\tilde{\varphi}_i \in \mathrm{H}^2(\Omega)$, $i = 1\ldots N$. Moreover, the coefficients $c_i$ in (21) are given by the following formula*

$$c_i = -(\pi)^{-1}(g, \zeta_i)_\Omega,$$

*where the $\zeta_i$ are defined in (19).*

Noticing that $\boldsymbol{x} \mapsto r_i^{\pi/\alpha_i} \sin(\pi\theta_i/\alpha_i)$ belongs to $\mathrm{H}^1(\Omega)\backslash\mathrm{H}^2(\Omega)$, we deduce the

**Corollary 2.2.** *Let $\varphi \in \mathrm{H}_0^1(\Omega)$ be such that $\Delta\varphi = g \in \mathrm{L}^2(\Omega)$. Then $\varphi \in \mathrm{H}^2(\Omega)$ if and only if $(g, \zeta_i)_\Omega = 0$ for $i = 1\ldots N$.*

We are now in position to start the investigation of $(\tilde{\mathscr{P}})$. To get an idea, let us consider $v$ a solution of $(\tilde{\mathscr{P}})$ with $f = 0$. According to Proposition 2.6, we have $\sigma\Delta v = \sum_{i=1}^N a_i \zeta_i$ where $a_1, \ldots, a_N$ are some constants. Since $v \in \mathrm{H}^2(\Omega)$, we must have $(\sum_{i=1}^N a_i \sigma^{-1}\zeta_i, \zeta_j)_\Omega = 0$ for $j = 1\ldots N$. This leads us to introduce the $N \times N$ matrix ($N$ is the number of reentrant corners)

$$\mathbb{M} := \begin{pmatrix} (\sigma^{-1}\zeta_1, \zeta_1)_\Omega & \cdots & (\sigma^{-1}\zeta_1, \zeta_N)_\Omega \\ \vdots & \ddots & \vdots \\ (\sigma^{-1}\zeta_N, \zeta_1)_\Omega & \cdots & (\sigma^{-1}\zeta_N, \zeta_N)_\Omega \end{pmatrix}. \tag{22}$$

Noticing that the functions $\zeta_i$, $i = 1\ldots N$, are real valued, we deduce that $\mathbb{M}$ is a symmetric matrix. We divide our study according to the dimension of the kernel of $\mathbb{M}$.

<div align="center">

⋆ ⋆ ⋆ ⋆ ⋆

Case where $\mathbb{M}$ is invertible

⋆ ⋆ ⋆ ⋆ ⋆

</div>

We assume here that the matrix $\mathbb{M}$ is invertible. Let us construct a dual basis of $(\zeta_i)_{i=1}^N$. This will be helpful for the sequel.

**Lemma 2.1.** *Assume that $\sigma \in \mathrm{L}^\infty(\Omega)$ verifies $\sigma^{-1} \in \mathrm{L}^\infty(\Omega)$. Assume that $\sigma$ is such that the matrix $\mathbb{M}$ is invertible. Then there exist $N$ functions $\lambda_i$, $i = 1 \ldots N$, belonging to $\mathrm{span}(\zeta_1, \ldots, \zeta_N)$, such that*
$$(\sigma^{-1}\lambda_i, \zeta_j)_\Omega = \delta_{ij}, \qquad \text{for } i, j \in \{1, \ldots, N\}.$$

*Proof.* Let $\mathbb{A}$ denote the inverse of $\mathbb{M}$. It is sufficient to take $\lambda_i := \sum_{k=1}^N \mathbb{A}_{ik}\zeta_k$. □

Then, we can prove the

**Proposition 2.8.** *Assume that the boundary $\partial\Omega$ is a polygon which has exactly $N$ reentrant corners. Assume that $\sigma \in \mathrm{L}^\infty(\Omega)$ verifies $\sigma^{-1} \in \mathrm{L}^\infty(\Omega)$. Moreover, assume that $\sigma$ is such that the matrix $\mathbb{M}$ defined in (22) is invertible. Then, the operator $\tilde{B} : \mathrm{H}_0^1(\Omega) \cap \mathrm{H}^2(\Omega) \to \mathrm{H}_0^1(\Omega) \cap \mathrm{H}^2(\Omega)$ defined in (12) is an isomorphism.*

*Proof.* We introduce the operator $\mathtt{T}$ such that, for all $v \in \mathrm{H}_0^1(\Omega) \cap \mathrm{H}^2(\Omega)$, $\mathtt{T}v \in \mathrm{H}_0^1(\Omega)$ is the function such that $\Delta(\mathtt{T}v) = \sigma^{-1}(\Delta v + \sum_{i=1}^N a_i\lambda_i)$ with, for $i = 1 \ldots N$, $a_i = -(\sigma^{-1}\Delta v, \zeta_i)_\Omega$. Since $(\sigma^{-1}(\Delta v + \sum_{i=1}^N a_i\lambda_i), \zeta_j)_\Omega = 0$ for $j = 1 \ldots N$, we have $\mathtt{T}v \in \mathrm{H}^2(\Omega)$ according to Corollary 2.2. Thus, $\mathtt{T}$ is a continuous operator of $\mathrm{H}_0^1(\Omega) \cap \mathrm{H}^2(\Omega)$. Then, one checks as in the proof of Proposition 2.4 the equality $\tilde{B} \circ \mathtt{T} = \mathrm{Id}$. This proves that $\tilde{B}$ and $\mathtt{T}$ are isomorphisms of $\mathrm{H}_0^1(\Omega) \cap \mathrm{H}^2(\Omega)$ with $\tilde{B}^{-1} = \mathtt{T}$. □

**Remark 2.3.** *The assumption "$\mathbb{M}$ invertible" is only an assumption on the values of $\sigma$ with respect to the geometry of the domain. Indeed, the singularities $\zeta_i$, $i = 1 \ldots N$, only depend on the geometry of $\Omega$.*

**Remark 2.4.** *When there is only one reentrant corner in $\partial\Omega$, the matrix $\mathbb{M}$, a scalar in this case, is invertible if and only if $(\sigma^{-1}\zeta_1, \zeta_1)_\Omega \neq 0$. Thus, we find back the result of Proposition 2.4.*

**Remark 2.5.** *When $\sigma \geq C > 0$, the matrix $\mathbb{M}$ is always invertible. Indeed, if $(a_1 \ldots a_N)^t \in \ker \mathbb{M}$, then the function $\tau := \sum_{i=1}^N a_i\zeta_i$ verifies $(\sigma^{-1}\tau, \zeta_j)_\Omega = 0$ for $j = 1 \ldots N$, and so $(\sigma^{-1}\tau, \tau)_\Omega = 0$. This implies $\tau = 0$. Since the family $(\zeta_i)_{i=1}^N$ is free, we deduce that $a_1 = \cdots = a_N = 0$. This is coherent with the fact that $\tilde{b}$ is coercive on $\mathrm{H}_0^1(\Omega) \cap \mathrm{H}^2(\Omega) \times \mathrm{H}_0^1(\Omega) \cap \mathrm{H}^2(\Omega)$ when $\sigma \geq C > 0$.*

<div align="center">

⋆ ⋆ ⋆ ⋆ ⋆

CASE WHERE $\mathbb{M}$ IS NOT INVERTIBLE

⋆ ⋆ ⋆ ⋆ ⋆

</div>

Now, we assume that $\sigma$ is such that the matrix $\mathbb{M}$ defined in (22) has a non trivial kernel of dimension $M > 0$. In order to study the properties of $(\tilde{\mathscr{P}})$ is such configurations, we need first to present some simple results of elementary linear algebra. Let us introduce $(\vec{\beta}_m)_{m=1}^M$ a basis of $\ker \mathbb{M}$, with, for $m = 1 \ldots M$, $\vec{\beta}_m = (\beta_{m1} \ldots \beta_{mN})^t$. Next, we define the functions

$$\beta_m := \sum_{i=1}^N \beta_{mi}\zeta_i \in \mathrm{L}^2(\Omega), \qquad m = 1 \ldots M. \tag{23}$$

**Lemma 2.2.** *The family $(\beta_m)_{m=1}^M$ constitutes a basis of the set*

$$\mathfrak{X} := \{\zeta \in span(\zeta_1, \ldots, \zeta_N) \,|\, (\sigma^{-1}\zeta, \zeta_j)_\Omega = 0, \text{ for } j = 1 \ldots N\}. \tag{24}$$

*Proof.* Let $\beta_m$ be defined as in (23). For $j = 1 \ldots N$, since $\vec{\beta}_m \in \ker \mathbb{M}$, we find $(\sigma^{-1}\beta_m, \zeta_j)_\Omega = \sum_{i=1}^N \beta_{mi}(\sigma^{-1}\zeta_i, \zeta_j)_\Omega = 0$ (recall that $\mathbb{M}$ is symmetric). We deduce that, for $m = 1 \ldots M$, $\beta_m \in \mathfrak{X}$. Now, we prove that $(\beta_m)_{m=1}^M$ forms a free family. Let $a_1, \ldots, a_M$ be some constants such that $\sum_{m=1}^M a_m \beta_m = 0$. Replacing $\beta_m$ by its expression (23), we obtain

$$0 = \sum_{m=1}^M a_m \Big( \sum_{i=1}^N \beta_{mi} \zeta_i \Big) = \sum_{i=1}^N \Big( \sum_{m=1}^M a_m \beta_{mi} \Big) \zeta_i = 0.$$

Since the family $(\zeta_i)_{i=1}^N$ forms a linearly independent set, it yields $\sum_{m=1}^M a_m \beta_{mi} = 0$ for $i = 1 \ldots N$. In other words, there holds $\sum_{m=1}^M a_m \vec{\beta}_m = 0$. But $(\vec{\beta}_m)_{m=1}^M$ is a basis of $\ker \mathbb{M}$. Therefore, we have $a_1 = \cdots = a_M = 0$. This proves that $(\beta_m)_{m=1}^M$ is free. Finally, if $\zeta = \sum_{i=1}^N a_i \zeta_i$ belongs to $\mathfrak{X}$, it is clear that $(a_1 \ldots a_N)^t$ is an element of $\ker \mathbb{M}$. This allows to conclude that $(\beta_m)_{m=1}^M$ is a basis of $\mathfrak{X}$. □

According to Lemma 2.2, we can find $N - M$ (linearly independent) functions $\gamma_1, \ldots, \gamma_{N-M}$ of $\mathrm{span}(\zeta_1, \ldots, \zeta_N)$ such that

$$\mathrm{span}(\zeta_1, \ldots, \zeta_N) = \mathrm{span}(\beta_1, \ldots, \beta_M) \oplus \mathrm{span}(\gamma_1, \ldots, \gamma_{N-M}). \tag{25}$$

In the next lemma, we construct a dual basis of $(\gamma_i)_{i=1}^{N-M}$.

**Lemma 2.3.** *Assume that $\sigma \in \mathrm{L}^\infty(\Omega)$ verifies $\sigma^{-1} \in \mathrm{L}^\infty(\Omega)$. Moreover, assume that $\sigma$ is such that the matrix $\mathbb{M}$ defined in (22) has a kernel of dimension $M > 0$. Then, there exist $N - M$ functions $\lambda_i$, $i = 1 \ldots N - M$, belonging to $\mathrm{span}(\gamma_1, \ldots, \gamma_{N-M})$, such that*

$$(\sigma^{-1}\lambda_i, \gamma_j)_\Omega = \delta_{ij}, \qquad \text{for } i, j \in \{1, \ldots, N-M\}.$$

*Proof.* Introduce the $(N - M) \times (N - M)$ matrix

$$\tilde{\mathbb{M}} := \begin{pmatrix} (\sigma^{-1}\gamma_1, \gamma_1)_\Omega & \cdots & (\sigma^{-1}\gamma_1, \gamma_{N-M})_\Omega \\ \vdots & \ddots & \vdots \\ (\sigma^{-1}\gamma_{N-M}, \gamma_1)_\Omega & \cdots & (\sigma^{-1}\gamma_{N-M}, \gamma_{N-M})_\Omega \end{pmatrix}. \tag{26}$$

Let us prove that this matrix is invertible. If $(a_1, \ldots, a_{N-M})^t$ is an element of $\ker \tilde{\mathbb{M}}$, a simple computation shows that the function $\tau = \sum_{i=1}^{N-M} a_i \gamma_i$ satisfies $(\sigma^{-1}\tau, \gamma_j)_\Omega = 0$, for $j = 1 \ldots N-M$. Now, since $\tau$ is a linear combination of the $\zeta_i$ and since the $\beta_j$ belong to the space $\mathfrak{X}$ defined in (24), we know that there also holds $(\sigma^{-1}\tau, \beta_j)_\Omega = 0$ for $j = 1 \ldots M$. It follows from the decomposition (25) that $(\sigma^{-1}\tau, \zeta_j)_\Omega = 0$ for $j = 1 \ldots N$. This proves that $\tau$ is an element of $\mathfrak{X} = \mathrm{span}(\beta_1, \ldots, \beta_M)$ (Lemma 2.2). Thus, $\tau \in (\mathrm{span}(\beta_1, \ldots, \beta_M) \cap \mathrm{span}(\gamma_1, \ldots, \gamma_{N-M}))$. We deduce that $\tau = \sum_{i=1}^{N-M} a_i \gamma_i = 0$. But the family $(\gamma_i)_{i=1}^{N-M}$ is free. Therefore, $a_1 = \cdots = a_{N-M} = 0$. This shows that $\tilde{\mathbb{M}}$ is invertible.
Let us denote $\tilde{\mathbb{A}}$ the inverse of $\tilde{\mathbb{M}}$. It just remains to define $\lambda_i := \sum_{k=1}^{N-M} \mathbb{A}_{ik} \gamma_k$. □

When the matrix $\mathbb{M}$ has a kernel of dimension $M$, a kernel and a cokernel of dimension $M$ appear for $(\tilde{\mathscr{P}})$. To describe $\ker \tilde{B}$ and $\mathrm{coker}\, \tilde{B}$, we introduce, for $m = 1 \ldots M$, the function $\psi_m \in \mathrm{H}_0^1(\Omega)$ verifying

$$\Delta \psi_m = \sigma^{-1}\beta_m, \tag{27}$$

where $\beta_m$ is defined in (23). Since $(\sigma^{-1}\beta_m, \zeta_j)_\Omega = 0$ for $j = 1 \ldots N$, we know according to Corollary 2.2 that $\psi_m$ belongs to $\mathrm{H}^2(\Omega)$. Moreover, for all $v' \in \mathrm{H}_0^1(\Omega) \cap \mathrm{H}^2(\Omega)$, there holds $(\sigma \Delta \psi_m, \Delta v')_\Omega = (\beta_m, \Delta v')_\Omega = 0$. As a consequence, $\psi_m$ constitutes an element of $\ker \tilde{B}$.

**Proposition 2.9.** *Assume that the boundary $\partial \Omega$ is a polygon which has exactly $N$ reentrant corners. Assume that $\sigma \in \mathrm{L}^\infty(\Omega)$ verifies $\sigma^{-1} \in \mathrm{L}^\infty(\Omega)$. Moreover, assume that $\sigma$ is such that the matrix $\mathbb{M}$ has a kernel of dimension $M > 0$. Then,*

- $\dim(\ker \tilde{B}) = M$ *with* $\ker \tilde{B} = \mathrm{span}(\psi_1, \ldots, \psi_M)$;
- $\dim(\mathrm{coker}\,\tilde{B}) = M$ *and for all* $f \in (\mathrm{H}_0^1(\Omega) \cap \mathrm{H}^2(\Omega))^*$, $(\tilde{\mathscr{P}})$ *has a solution if and only if* $\langle f, \psi_m \rangle_\Omega = 0$ *for* $m = 1\ldots M$.

*In this statement, the matrix* $\mathbb{M}$ *and the functions* $\psi_m$ *are respectively defined in* (22) *and* (27).

*Proof.* ⋆ KERNEL. If $v$ belongs to $\ker \tilde{B}$ then, according to Proposition 2.6, there holds $\sigma \Delta v = \sum_{i=1}^N a_i \zeta_i$ where the $a_1, \ldots, a_N$ are some constants. Since $v \in \mathrm{H}^2(\Omega)$, according to Corollary 2.2, we have $(\sum_{i=1}^N a_i \sigma^{-1} \zeta_i, \zeta_j)_\Omega = 0$ for $j = 1 \ldots N$. This implies $(a_1 \ldots a_N)^t \in \ker \mathbb{M} = \mathrm{span}(\vec{\beta}_1, \ldots, \vec{\beta}_M)$. In other words, there exist constants $c_m$, $m = 1 \ldots M$, such that $(a_1 \ldots a_N)^t = \sum_{m=1}^M c_m \vec{\beta}_m$. It yields

$$\sigma \Delta v = \sum_{i=1}^N a_i \zeta_i = \sum_{i=1}^N \big(\sum_{m=1}^M c_m \beta_{mi}\big) \zeta_i = \sum_{m=1}^M c_m \big(\sum_{i=1}^N \beta_{mi} \zeta_i\big) = \sum_{m=1}^M c_m \beta_m.$$

As a consequence, we can write $v = \sum_{m=1}^M c_m \psi_m$. This proves the relation $\ker \tilde{B} \subset \mathrm{span}(\psi_1, \ldots, \psi_M)$. Since the reverse inclusion is true, we obtain $\ker \tilde{B} = \mathrm{span}(\psi_1, \ldots, \psi_M)$. Using (27), since $(\beta_m)_{m=1}^M$ is free, we deduce that $(\psi_m)_{m=1}^M$ is free. Therefore, $\ker \tilde{B} = \mathrm{span}(\psi_1, \ldots, \psi_M)$ is a space of dimension $M$.

⋆ COKERNEL. Consider $f \in (\mathrm{H}_0^1(\Omega) \cap \mathrm{H}^2(\Omega))^*$ such that $\langle f, \psi_m \rangle_\Omega = 0$ for $m = 1 \ldots M$. Let us introduce $p_0$ an element of $\mathrm{L}^2(\Omega)$ such that $(p_0, \Delta v')_\Omega = \langle f, v' \rangle_\Omega$ for all $v' \in \mathrm{H}_0^1(\Omega) \cap \mathrm{H}^2(\Omega)$. The existence of such a function $p_0$ is guaranteed by Proposition 2.6. Next, define $p = p_0 - \sum_{i=1}^{N-M} a_i \lambda_i$ with, for $i = 1 \ldots N - M$, $a_i := (\sigma^{-1} p_0, \gamma_i)_\Omega$. Here, the $\lambda_i$ are the functions defined by Lemma 2.3. By construction of the $a_i$, we have $(\sigma^{-1} p, \gamma_i)_\Omega = 0$ for $i = 1 \ldots N - M$. But we have also $(\sigma^{-1} p, \beta_m)_\Omega = 0$ for $m = 1 \ldots N$. Indeed, since $\Delta \psi_m = \sigma^{-1} \beta_m$, we can write $(\sigma^{-1} p, \beta_m)_\Omega = (p, \sigma^{-1} \beta_m)_\Omega = (p, \Delta \psi_m)_\Omega = \langle f, \psi_m \rangle_\Omega = 0$. Thus, the function $p$ satisfies the compatibility conditions $(\sigma^{-1} p, \zeta_i)_\Omega = 0$ for $i = 1 \ldots N$. As a consequence, by virtue of Corollary 2.2, the function $v \in \mathrm{H}_0^1(\Omega)$ verifying $\Delta v = \sigma^{-1} p$ is in $\mathrm{H}^2(\Omega)$. Moreover, for all $v' \in \mathrm{H}_0^1(\Omega) \cap \mathrm{H}^2(\Omega)$, there holds

$$(\sigma \Delta v, \Delta v')_\Omega = (p, \Delta v')_\Omega = (p_0 - \sum_{i=1}^{N-M} a_i \lambda_i, \Delta v')_\Omega = (p_0, \Delta v')_\Omega = \langle f, v' \rangle_\Omega.$$

Therefore, $v$ constitutes a solution to $(\tilde{\mathscr{P}})$. Now, let us consider $f \in (\mathrm{H}_0^1(\Omega) \cap \mathrm{H}^2(\Omega))^*$ such that $\langle f, \psi_m \rangle_\Omega \neq 0$ for some $m \in \{1, \ldots, M\}$. Let us assume that there exists a solution $v$ to $(\tilde{\mathscr{P}})$. Then, by Proposition 2.6, we have $\sigma \Delta v = p_0 + \sum_{j=1}^N a_j \zeta_j$. Here, $p_0 \in \mathrm{L}^2(\Omega)$ satisfies $(p_0, \Delta v')_\Omega = \langle f, v' \rangle_\Omega$ for all $v' \in \mathrm{H}_0^1(\Omega) \cap \mathrm{H}^2(\Omega)$ and the $a_1, \ldots, a_N$ are some constants. This imposes $(\sigma^{-1} p_0, \beta_m)_\Omega = 0$ for $m = 1 \ldots M$. But there holds $(\sigma^{-1} p_0, \beta_m)_\Omega = (p_0, \sigma^{-1} \beta_m)_\Omega = (p_0, \Delta \psi_m)_\Omega = \langle f, \psi_m \rangle_\Omega$. Thus, we obtain an absurdity. This ends to prove that there exists a solution to $(\tilde{\mathscr{P}})$ if and only if $\langle f, \psi_m \rangle_\Omega = 0$ for $m = 1 \ldots M$. □

⋄ EXAMPLE. Consider the open set $\Omega$ of Figure 3 which presents the particularity to be symmetric with respect to the axis $(Ox)$. Notice also that the vertices of the reentrant corners are located on $(Ox)$. For this configuration, it is easy to prove that the functions $\zeta_1$ and $\zeta_2$ defined in (19) are symmetric with respect to $(Ox)$. Indeed, for $i = 1, 2$, the function $\hat{\zeta}_i : (x, y) \mapsto \zeta_i(x, -y)$ verifies Problem (20) with $f = 0$. But Proposition 2.6 indicates that $(\zeta_i)_{i=1}^2$ is a basis of the set of functions satisfying Problem (20) with $f = 0$. Using the behaviour of $\hat{\zeta}_i$ at $O_i$, this allows us to conclude that $\hat{\zeta}_i = \zeta_i$ for $i = 1, 2$.

Therefore, when $\sigma$ is skewsymmetric with respect to the axis $(Ox)$, *i.e.* when $\sigma(x, y) = -\sigma(x, -y)$ a.e. in $\Omega$, the matrix $\mathbb{M}$ is equal to the null matrix. In this situation, $(\tilde{\mathscr{P}})$ has a kernel and a cokernel which are both of dimension two.

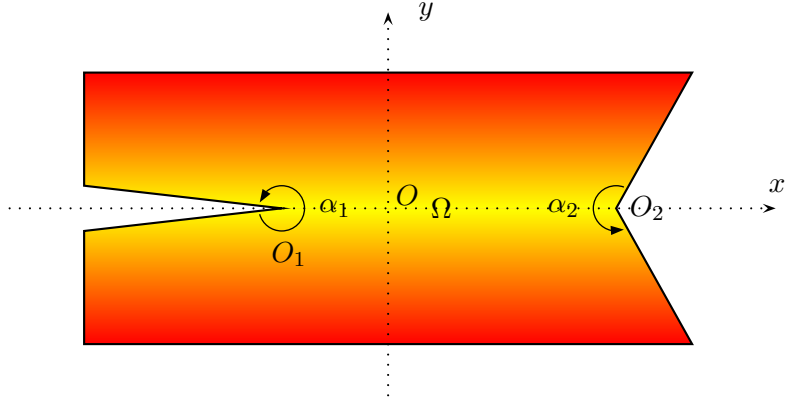

Figure 3: An example of a domain with a polygonal boundary which is symmetric with respect to the axis $(Ox)$. Here, $\Omega$ has two reentrant corners – $N = 2$.

## 2.3 Problem set in $\mathrm{H}_0^1(\Delta)$

Up to now in this section, we have imposed the boundary condition $\sigma \Delta v = 0$ on $\partial \Omega$ in a weak way working with the integral identity "Find $v \in \mathrm{H}_0^1(\Omega) \cap \mathrm{H}^2(\Omega)$ such that $(\sigma \Delta v, \Delta v')_\Omega = \langle f, v' \rangle_\Omega$ for all $v' \in \mathrm{H}_0^1(\Omega) \cap \mathrm{H}^2(\Omega)$". Notice that this variational formulation has also a sense, for a source term smooth enough, when the functions $v$, $v'$ are chosen in the space $\mathrm{H}_0^1(\Delta) := \{\varphi \in \mathrm{H}_0^1(\Omega) \,|\, \Delta \varphi \in \mathrm{L}^2(\Omega)\}$. Let $\Omega$ be a domain (with a Lipschitz boundary) of $\mathbb{R}^d$, $d \geq 1$. In this paragraph, we wish to study the problem,

$$(\mathscr{P}^\sharp) \left|\begin{array}{l} \text{Find } v \in \mathrm{H}_0^1(\Delta) \text{ such that:} \\ (\sigma \Delta v, \Delta v')_\Omega = \langle f, v' \rangle_\Omega, \qquad \forall v' \in \mathrm{H}_0^1(\Delta). \end{array}\right. \tag{28}$$

Here, $f$ is a given source term of $\mathrm{H}_0^1(\Delta)^*$, the topological dual space of $\mathrm{H}_0^1(\Delta)$. According to the Lax-Milgram theorem, we have $\|\varphi\|_{\mathrm{H}_0^1(\Omega)} \leq C \|\Delta \varphi\|_\Omega$ for all $\varphi \in \mathrm{H}_0^1(\Delta)$. Therefore, $(v, v') \mapsto (v, v')_{\mathrm{H}_0^2(\Omega)} = (\Delta v, \Delta v')_\Omega$ defines an inner product on $\mathrm{H}_0^1(\Delta)$ and the associated norm is equivalent to the natural norm $v \mapsto (\|v\|_{\mathrm{H}_0^1(\Omega)}^2 + \|\Delta v\|_\Omega^2)^{1/2}$. In order to study the properties of $(\mathscr{P}^\sharp)$, we introduce with the Riesz representation theorem, the bounded operator $B^\sharp : \mathrm{H}_0^1(\Delta) \to \mathrm{H}_0^1(\Delta)$ such that

$$(B^\sharp v, v')_{\mathrm{H}_0^2(\Omega)} = (\sigma \Delta v, \Delta v'), \qquad \forall (v, v') \in \mathrm{H}_0^1(\Delta) \times \mathrm{H}_0^1(\Delta). \tag{29}$$

There holds the

**Proposition 2.10.** *Assume that $\sigma \in \mathrm{L}^\infty(\Omega)$ is such that $\sigma^{-1} \in \mathrm{L}^\infty(\Omega)$. Then the operator $B^\sharp : \mathrm{H}_0^1(\Delta) \to \mathrm{H}_0^1(\Delta)$ defined in (29) is an isomorphism.*

*Proof.* Let us introduce the operator $\mathtt{T}$ such that, for $v \in \mathrm{H}_0^1(\Delta)$, $\mathtt{T}v \in \mathrm{H}_0^1(\Delta)$ is the function satisfying $\Delta(\mathtt{T}v) = \sigma^{-1} \Delta v$. It is clear that $B^\sharp$ is a continuous operator. For all $(v, v') \in \mathrm{H}_0^1(\Delta) \times \mathrm{H}_0^1(\Delta)$, we have

$$(B^\sharp(\mathtt{T}v), v')_{\mathrm{H}_0^2(\Omega)} = (\sigma \Delta(\mathtt{T}v), \Delta v')_\Omega = (\Delta v, \Delta v')_\Omega.$$

This proves that $B^\sharp \circ \mathtt{T} = \mathrm{Id}$. Since $B^\sharp$ is selfadjoint, we deduce that $B^\sharp$ is an isomorphism. Its inverse is equal to $\mathtt{T}$. $\square$

Let us study the smoothness of the solution $v^\sharp$ of Problem $(\mathscr{P}^\sharp)$ defined in (28). More precisely, our goal is to compare $v^\sharp$ with $\tilde{v}$, the solution of Problem $(\tilde{\mathscr{P}})$ defined in (10), when the latter is well-posed. Since $(\mathrm{H}_0^1(\Omega) \cap \mathrm{H}^2(\Omega)) \subset \mathrm{H}_0^1(\Delta)$, there holds $\mathrm{H}_0^1(\Delta)^* \subset (\mathrm{H}_0^1(\Omega) \cap \mathrm{H}^2(\Omega))^*$. Therefore, the source term $f$ can be chosen in $\mathrm{H}_0^1(\Delta)^*$.

As soon as the domain $\Omega$ is such that the operator $\Delta : \mathrm{H}_0^1(\Omega) \cap \mathrm{H}^2(\Omega) \to \mathrm{L}^2(\Omega)$ is an isomorphism, we have $\mathrm{H}_0^1(\Delta) = \mathrm{H}_0^1(\Omega) \cap \mathrm{H}^2(\Omega)$. In this case, for example when $\Omega \subset \mathbb{R}^d$ is convex or of class $\mathscr{C}^2$, there holds $v^\sharp = \tilde{v}$.

Now, let us study the situation where $\mathrm{H}_0^1(\Delta) \neq \mathrm{H}_0^1(\Omega) \cap \mathrm{H}^2(\Omega)$. Let us work in dimension 2, for a domain $\Omega$ which has, to set ideas, exactly one reentrant corner located at $O$. We reintroduce the function $\zeta$ defined in (14) verifying $\zeta \in \mathrm{L}^2(\Omega) \backslash \mathrm{H}^1(\Omega)$, $\Delta \zeta = 0$ a.e. in $\Omega$ and $\zeta = 0$ a.e. on $\partial \Omega$. As in (18), we define $\psi \in \mathrm{H}_0^1(\Omega)$ the function such that $\Delta \psi = \sigma^{-1} \zeta$. Since $\zeta \in \mathrm{L}^2(\Omega)$, we have $\psi \in \mathrm{H}_0^1(\Delta)$ and so $(\sigma \Delta v^\sharp, \Delta \psi)_\Omega = \langle f, \psi \rangle_\Omega$. This can also be written $(\Delta v^\sharp, \zeta)_\Omega = \langle f, \psi \rangle_\Omega$. Thus, by virtue of Corollary 2.1, we have $v^\sharp \in \mathrm{H}_0^1(\Omega) \cap \mathrm{H}^2(\Omega)$ if and only if $\langle f, \psi \rangle_\Omega = 0$. Let us distinguish two cases.

- If $(\sigma^{-1} \zeta, \zeta)_\Omega \neq 0$, according to Proposition 2.4, the operator $\tilde{B}$ defined in (12) is an isomorphism. Therefore, the solution $\tilde{v}$ of $(\tilde{\mathscr{P}})$ is defined uniquely.
  - If $\langle f, \psi \rangle_\Omega = 0$, then $v^\sharp$ verifies the same problem as $\tilde{v}$. We deduce that $v^\sharp = \tilde{v}$ in this configuration.
  - If $\langle f, \psi \rangle_\Omega \neq 0$, then $v^\sharp \notin \mathrm{H}^2(\Omega)$ and so $v^\sharp \neq \tilde{v}$. More precisely, since $(\sigma \Delta(v^\sharp - \tilde{v}), \Delta v')_\Omega = 0$ for all $v' \in \mathrm{H}_0^1(\Omega) \cap \mathrm{H}^2(\Omega)$, we deduce that $\Delta(v^\sharp - \tilde{v}) = a\, \sigma^{-1} \zeta$, where $a$ is constant. Multiplying by $\zeta$, integrating by parts on $\Omega$ and using that $\tilde{v} \in \mathrm{H}_0^1(\Omega) \cap \mathrm{H}^2(\Omega)$, we deduce that $a = (\Delta v^\sharp, \zeta)_\Omega / (\sigma^{-1} \zeta, \zeta)_\Omega = \langle f, \psi \rangle_\Omega / (\sigma^{-1} \zeta, \zeta)_\Omega$. Thus, in this configuration, we have
  $$v^\sharp - \tilde{v} = \frac{\langle f, \psi \rangle_\Omega}{(\sigma^{-1} \zeta, \zeta)_\Omega} \psi.$$

- If $(\sigma^{-1} \zeta, \zeta)_\Omega = 0$, then, according to Proposition 2.5, $(\tilde{\mathscr{P}})$ has a solution if and only if $\langle f, \psi \rangle_\Omega = 0$. Let us assume that this holds true. In this case, according to the above discussion, the solution $v^\sharp \in \mathrm{H}_0^1(\Delta)$ of $(\mathscr{P}^\sharp)$ is actually in $\mathrm{H}_0^1(\Omega) \cap \mathrm{H}^2(\Omega)$. Now, let us explain why $\psi$ belongs to $\ker \tilde{B}$ but not to $\ker B^\sharp$. A function $v \in \mathrm{H}_0^1(\Omega) \cap \mathrm{H}^2(\Omega)$ belongs to $\ker \tilde{B}$ if and only if it satisfies $(\sigma \Delta v, \Delta v')_\Omega = 0$ for all $v' \in \mathrm{H}_0^1(\Omega) \cap \mathrm{H}^2(\Omega)$. Since $(\sigma \Delta \psi, \Delta v')_\Omega = (\zeta, \Delta v')_\Omega$, we have $\psi \in \ker \tilde{B}$. On the other hand, $\psi \in \ker B^\sharp$ if and only if we have $(\sigma \Delta \psi, \Delta v')_\Omega = 0$ for all $v' \in \mathrm{H}_0^1(\Delta) \supset (\mathrm{H}_0^1(\Omega) \cap \mathrm{H}^2(\Omega))$. Testing against $v'$ such that $\Delta v' = \zeta$, we see that $\psi$ is not an element of $\ker B^\sharp$. Of course, such a $v'$ does not belong to $\mathrm{H}_0^1(\Omega) \cap \mathrm{H}^2(\Omega)$ because $(\zeta, \zeta)_\Omega \neq 0$.

## 3 Bilaplacian with Dirichlet boundary condition

In this section, we come back to the study of the operator $B : \mathrm{H}_0^2(\Omega) \to \mathrm{H}_0^2(\Omega)$ introduced in (7). First, we provide a sufficient criterion to ensure that $B$ is Fredholm of index zero. Then, we exhibit situations where $B$ is not of Fredholm type.

### 3.1 Configurations where $\sigma$ has a constant sign on the boundary

#### 3.1.1 Fredholm property

In this paragraph, $\Omega$ is a domain of $\mathbb{R}^d$, with $d \geq 1$. We prove that $B$ is Fredholm of index zero when $\sigma$ satisfies the following condition.

$(\mathscr{H}_\sigma)$ We assume that $\sigma \in \mathrm{L}^\infty(\Omega)$ is such that $\sigma^{-1} \in \mathrm{L}^\infty(\Omega)$. Moreover, we assume that $\sigma(\boldsymbol{x}) \geq C_1 > 0$ a.e. in $\Omega \backslash \overline{\mathcal{O}}$ or $\sigma(\boldsymbol{x}) \leq C_2 < 0$ a.e. in $\Omega \backslash \overline{\mathcal{O}}$, where $C_1$, $C_2$ are two constants and where $\mathcal{O}$ is an open set such that $\overline{\mathcal{O}} \subset \Omega$.

In other words, we assume that there exists a neighbourhood of $\partial \Omega$ where $\sigma \geq C_1 > 0$ or $\sigma \leq C_2 < 0$. However, outside this neighbourhood, $\sigma$ can change sign. To prove that $B$ is a Fredholm operator of index zero, we build a right parametrix for $B$, *i.e.* we build a bounded operator $\mathtt{T}$ such that $B\circ\mathtt{T} = \mathcal{I} + \mathcal{K}$ where $\mathcal{I} : \mathrm{H}_0^2(\Omega) \to \mathrm{H}_0^2(\Omega)$ is an isomorphism and $\mathcal{K} : \mathrm{H}_0^2(\Omega) \to \mathrm{H}_0^2(\Omega)$ is a compact operator.

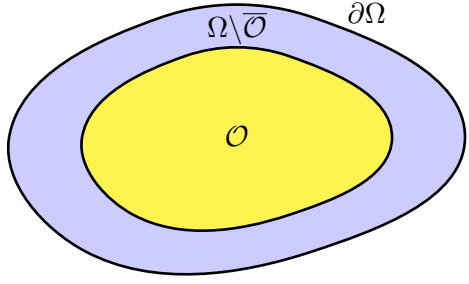

Figure 4: The parameter $\sigma$ is assumed to be uniformly positive or uniformly negative in the region $\Omega\setminus\overline{\mathcal{O}}$.

**Theorem 3.1.** *Assume that $\sigma$ satisfies condition $(\mathscr{H}_\sigma)$. Then the operator $B : \mathrm{H}_0^2(\Omega) \to \mathrm{H}_0^2(\Omega)$ verifying $(Bv, v')_{\mathrm{H}_0^2(\Omega)} = (\sigma\Delta v, \Delta v')_\Omega$, for all $(v, v') \in \mathrm{H}_0^2(\Omega) \times \mathrm{H}_0^2(\Omega)$, is Fredholm of index zero.*

**Remark 3.1.** *Using the approach presented here, we can prove that for all $m \in \mathbb{N}^* := \{1, 2, 3, \dots\}$ the operator $\Delta^m(\sigma\Delta^m\cdot) : \mathrm{H}_0^{2m}(\Omega) \to \mathrm{H}^{-2m}(\Omega)$ is Fredholm of index zero when $\sigma$ satisfies condition $(\mathscr{H}_\sigma)$. Here, $\mathrm{H}_0^{2m}(\Omega)$ denotes the closure of $\mathscr{C}_0^\infty(\Omega)$ for the $\mathrm{H}^{2m}$-norm.*

*Proof.* Let us give the proof in the case where $\sigma \geq C_1 > 0$ in a neighbourhood of the boundary $\partial\Omega$. The configuration where $\sigma \leq C_2 < 0$ in a neighbourhood of $\partial\Omega$ can be deduced easily working with the operator $-B$. Let us introduce $\zeta \in \mathscr{C}_0^\infty(\Omega, [0; 1])$ a cut-off function equal to 1 in $\mathcal{O}$. Notice that $1 - \zeta$ is an element of $\mathscr{C}^\infty(\overline{\Omega}, [0; 1])$ which is equal to 1 in a neighbourhood of $\partial\Omega$. Now, let us consider $v$ an element of $\mathrm{H}_0^2(\Omega)$. The function $(1-\zeta)v$ belongs to $\mathrm{H}_0^2(\Omega)$ and, by definition, we have, for all $v' \in \mathrm{H}_0^2(\Omega)$,

$$(\sigma\Delta((1-\zeta)v), \Delta v')_\Omega = b((1-\zeta)v, v') = (B((1-\zeta)v), v')_{\mathrm{H}_0^2(\Omega)}.$$

This allows to write, expanding $\Delta((1-\zeta)v)$,

$$((1-\zeta)\sigma\Delta v, \Delta v')_\Omega = (B((1-\zeta)v), v')_{\mathrm{H}_0^2(\Omega)} + (\sigma(2\nabla v \cdot \nabla\zeta + v\Delta\zeta), \Delta v')_\Omega. \qquad (30)$$

On the support of $\zeta$, we must proceed slightly differently because we allow $\sigma$ to change sign. Let us denote $\psi$ the unique element of $\mathrm{H}_0^1(\Omega)$ such that $\Delta\psi = \sigma^{-1}\Delta v \in \mathrm{L}^2(\Omega)$. The classical results of interior regularity (see for example [24, theorem 2.1.3]) indicate that, for all $\chi \in \mathscr{C}_0^\infty(\Omega)$, $\chi\psi \in \mathrm{H}_0^2(\Omega)$ with the estimate $\|\chi\psi\|_{\mathrm{H}_0^2(\Omega)} \leq C\|\sigma^{-1}\Delta v\|_\Omega \leq C\|v\|_{\mathrm{H}_0^2(\Omega)}$. In particular, the function $\zeta\psi$ belongs to $\mathrm{H}_0^2(\Omega)$ and depends continuously on $v$. Since $(\sigma\Delta(\zeta\psi), \Delta v')_\Omega = (B(\zeta\psi), v')_{\mathrm{H}_0^2(\Omega)}$, we obtain, expanding $\Delta(\zeta\psi)$,

$$(\zeta\Delta v, \Delta v')_\Omega = (\sigma\zeta\Delta\psi, \Delta v')_\Omega = (B(\zeta\psi), v')_{\mathrm{H}_0^2(\Omega)} - (\sigma(2\nabla\psi \cdot \nabla\zeta + \psi\Delta\zeta), \Delta v')_\Omega. \qquad (31)$$

Let us define the operator $\mathtt{T} : \mathrm{H}_0^2(\Omega) \to \mathrm{H}_0^2(\Omega)$ such that $\mathtt{T}v = \zeta\psi + (1-\zeta)v$ for all $v \in \mathrm{H}_0^2(\Omega)$. With the Riesz representation, we introduce the operators $\mathcal{I} : \mathrm{H}_0^2(\Omega) \to \mathrm{H}_0^2(\Omega)$ and $\mathcal{K} : \mathrm{H}_0^2(\Omega) \to \mathrm{H}_0^2(\Omega)$ such that, for all $(v, v') \in \mathrm{H}_0^2(\Omega) \times \mathrm{H}_0^2(\Omega)$,

$$\begin{aligned}(\mathcal{I}v, v')_{\mathrm{H}_0^2(\Omega)} &= ((\zeta + (1-\zeta)\sigma)\Delta v, \Delta v')_\Omega \\ (\mathcal{K}v, v')_{\mathrm{H}_0^2(\Omega)} &= (\sigma(2\nabla(\psi-v)\cdot\nabla\zeta + (\psi-v)\Delta\zeta), \Delta v')_\Omega.\end{aligned} \qquad (32)$$

With these definitions, we have the relation $B\circ\mathtt{T} = \mathcal{I} + \mathcal{K}$. One verifies easily that $\mathcal{I} : \mathrm{H}_0^2(\Omega) \to \mathrm{H}_0^2(\Omega)$ is an isomorphism because the sesquilinear form $(v, v') \mapsto ((\zeta + (1-\zeta)\sigma)\Delta v, \Delta v')_\Omega$ is coercive on $H_0^2(\Omega) \times H_0^2(\Omega)$. Lemma 3.1 hereafter indicates that $\mathcal{K} : H_0^2(\Omega) \to H_0^2(\Omega)$ is compact. Therefore, the operator T constitutes a right parametrix for $B$. Since $B$ is selfadjoint, we deduce that it is a Fredholm operator of index zero. □

**Lemma 3.1.** *The operator $\mathcal{K} : H_0^2(\Omega) \to H_0^2(\Omega)$ defined in (32) is compact.*

*Proof.* Let us consider $(v'_m)_m$ a bounded sequence of elements of $H_0^2(\Omega)$. Let us prove that we can extract a subsequence of $(v'_m)_m$, still denoted $(v'_m)_m$, such that $(\mathcal{K} v'_m)_m$ converges in $H_0^2(\Omega)$. According to the definition (32) of $\mathcal{K}$, we can write

$$\|\mathcal{K} v'_m\|^2_{H_0^2(\Omega)} \leq C \, \|v'_m\|_{H_0^2(\Omega)} \left( \|v'_m\|_{H^1(\Omega)} + \|\psi_m\|_{H^1(\operatorname{supp}\zeta)} \right).$$

In the above equation, "$\operatorname{supp}\zeta$" designates the support of $\zeta$. Let us remind that by virtue of the result of interior regularity, there holds the estimate $\|\psi_m\|_{H^2(\operatorname{supp}\zeta)} \leq C \, \|v'_m\|_{H_0^2(\Omega)}$. Since the embedding of $H^2(\Omega)$ (resp. $H^2(\operatorname{supp}\zeta)$) in $H^1(\Omega)$ (resp. $H^1(\operatorname{supp}\zeta)$) is compact, we can extract a subsequence of $(v'_m)_m$, still denoted $(v'_m)_m$, such that $(v'_m)_m$ and $(\psi_m)_m$ converge respectively in $H^1(\Omega)$ and in $H^1(\operatorname{supp}\zeta)$ strongly. Let us define $v'_{mp} := v'_m - v'_p$, $\psi_{mp} := \psi_m - \psi_p$ for all $m$, $p \in \mathbb{N}$. Using the estimate

$$\|\mathcal{K} v'_{mp}\|^2_{H_0^2(\Omega)} \leq C \, \|v'_{mp}\|_{H_0^2(\Omega)} \left( \|v'_{mp}\|_{H^1(\Omega)} + \|\psi_{mp}\|_{H^1(\operatorname{supp}\zeta)} \right),$$

we deduce that $(\mathcal{K} v'_m)_m$ is a Cauchy sequence of $H_0^2(\Omega)$. As a consequence, $(\mathcal{K} v'_m)_m$ converges. □

Recall that we denoted $A_k$ the operator associated with the original interior transmission problem (see the definition in (5)). For all $k \in \mathbb{C}$, $A_k - B$ is a compact operator. Since the index of an operator is stable under compact perturbations (see for example [46, theorem 12.8]), we have the

**Corollary 3.1.** *Assume that $\sigma$ satisfies condition $(\mathscr{H}_\sigma)$. Then, for all $k \in \mathbb{C}$, the operator $A_k$ defined in (5) is Fredholm of index zero.*

**Remark 3.2.** *In a quite surprising way, this result indicates that Fredholmness for the operator $\Delta(\sigma\Delta\cdot) : H_0^2(\Omega) \to H^{-2}(\Omega)$ does not depend on the changes of sign of $\sigma$ which occur inside the domain $\Omega$. This property is not true for the operator $\operatorname{div}(\sigma\nabla\cdot) : H_0^1(\Omega) \to H^{-1}(\Omega)$ (see [19, 7, 2]).*

### 3.1.2 Study of the injectivity in 1D

In this paragraph, we wish to know whether the result we just obtained is optimal or not. More precisely, we proved that the operator $B$ is Fredholm of index zero when $\sigma$ remains positive or negative in a neighbourhood of $\partial\Omega$. As for the operator $\tilde{B}$ in convex or smooth domains, we may have a stronger property. Maybe $B$ is an isomorphism of $H_0^2(\Omega)$ as soon as $\sigma$ remains positive or negative in a neighbourhood of $\partial\Omega$. We will see on 1D examples for which we can carry explicit computations that this is not true: even for convex domains, $B$ can have a non trivial kernel.

⋄ EXAMPLE 1. Let us define the domains $\Omega = (a; b)$, $\Omega_1 = (a; 0)$, $\Omega_2 = (0; b)$, with $a < 0$ and $b > 0$. We introduce the function $\sigma$ such that $\sigma = \sigma_1$ in $\Omega_1$, $\sigma = \sigma_2$ in $\Omega_2$. Here, $\sigma_1 > 0$ and $\sigma_2 < 0$ are some constants. Using the proof of Theorem 3.1, one shows that in this configuration, the operator $B$ is Fredholm of index zero. Therefore, to know whether or not $B$ is an isomorphism, it is sufficient to study $\ker B$. In the sequel, if $v$ is an element of $H_0^2(\Omega)$, we denote $v_i := v|_{\Omega_i}$ for $i = 1, 2$. Classically, one proves that if $v \in H_0^2(\Omega)$ satisfies $(\mathscr{P})$ with $f = 0$, then $(v_1, v_2) \in H^2(\Omega_1) \times H^2(\Omega_2)$ verifies the transmission problem

$$\left|\begin{array}{lcl} \sigma_1 \Delta\Delta v_1 & = & 0 \quad \text{in } \Omega_1 \\ \sigma_1 \Delta\Delta v_2 & = & 0 \quad \text{in } \Omega_2 \\ v_1(a) = v_2(b) = v_1^{(1)}(a) = v_2^{(1)}(b) & = & 0 \\ v_1(0) - v_2(0) = v_1^{(1)}(0) - v_2^{(1)}(0) & = & 0 \\ \sigma_1 v_1^{(2)}(0) - \sigma_2 v_2^{(2)}(0) = \sigma_1 v_1^{(3)}(0) - \sigma_2 v_2^{(3)}(0) & = & 0 \end{array}\right. \quad (33)$$

In the above problem, $v^{(k)}(x)$ designates the $k$th derivative of $v$ at point $x$. Of course, if $(v_1, v_2) \in \mathrm{H}^2(\Omega_1) \times \mathrm{H}^2(\Omega_2)$ satisfies (33), then the function $v$ such that $v|_{\Omega_i} = v_i$, for $i = 1, 2$, is an element of $\ker B$. Using the first three lines of (33), we can write

$$v_1(x) = A_1(x-a)^3 + B_1(x-a)^2 \text{ for } x \in \Omega_1 \quad \text{and} \quad v_2(x) = A_2(x-b)^3 + B_2(x-b)^2 \text{ for } x \in \Omega_2.$$

The last two lines of (33) impose:

$$\begin{array}{rclcrcl}
-a^3 A_1 + a^2 B_1 & = & -b^3 A_2 + b^2 B_2 \,; & & 3a^2 A_1 - 2a B_1 & = & 3b^2 A_2 - 2b B_2 \,; \\
\sigma_1(-6a A_1 + 2 B_1) & = & \sigma_2(-6b A_2 + 2 B_2) \,; & & 6\sigma_1 A_1 & = & 6\sigma_2 A_2.
\end{array}$$

We deduce that $\ker B$ is non trivial if and only if the contrast $\kappa_\sigma := \sigma_2/\sigma_1$ satisfies

$$\kappa_\sigma^2 + \left(-4(b/a) + 6(b/a)^2 - 4(b/a)^3\right) \kappa_\sigma + (b/a)^4 = 0.$$

We can check that the discriminant of this polynomial remains positive for $(b/a) \in \mathbb{R}_-^* := (-\infty; 0)$. Thus, for all $(b/a) \in \mathbb{R}_-^*$, there exist two values of the contrast

$$\kappa_\sigma = \left(2 - 3(b/a) + 2(b/a)^2 \pm 2|(b/a) - 1|\sqrt{((b/a)^2 - (b/a) + 1)}\right)(b/a)$$

for which there exists a non trivial solution to (33). By a straightforward computation, one proves that these two roots are strictly negative for $(b/a) \in \mathbb{R}_-^*$. This is rather reassuring because $(v, v') \mapsto (\sigma \Delta v, \Delta v')_\Omega$ is coercive on $\mathrm{H}_0^2(\Omega) \times \mathrm{H}_0^2(\Omega)$ when $\kappa_\sigma > 0$. In the case where the domain is symmetric with respect to $x = 0$, that is, in the case where $b = -a$, $B$ is not injective for $\kappa_\sigma = -4 \pm \sqrt{3}$. Thus, these computations show that, according to the values of $\sigma$, the Fredholm operator $B : \mathrm{H}_0^2(\Omega) \to \mathrm{H}_0^2(\Omega)$ is *not always injective*.

$\diamond$ EXAMPLE 2. Let us look at a configuration where $\sigma$ has a constant sign in a neighbourhood of $\partial \Omega$. Define the open sets $\Omega = (-1; 1)$, $\Omega_1 = (-1; -\delta) \cup (\delta; 1)$, $\Omega_2 = (-\delta; \delta)$, with $0 < \delta < 1$. We introduce the function $\sigma$ such that $\sigma = \sigma_1$ in $\Omega_1$, $\sigma = \sigma_2$ in $\Omega_2$. Again, $\sigma_1 > 0$ and $\sigma_2 < 0$ are some constants. According to Theorem 3.1, for all contrast $\kappa_\sigma \in \mathbb{R}_-^*$, the operator $B : \mathrm{H}_0^2(\Omega) \to \mathrm{H}_0^2(\Omega)$ is Fredholm of index zero. Proceeding as for Example 1, we find that it is an isomorphism if and only if

$$\kappa_\sigma \notin \{\delta^3/(\delta^3 - 1), \delta/(\delta - 1)\}.$$

Again, this proves that the result of the previous paragraph is not under-optimal in the sense that $B$ is not always an isomorphism of $\mathrm{H}_0^2(\Omega)$.

To apply the analytic Fredholm theorem to prove that the set of transmission eigenvalues is discrete, we need to find some $k \in \mathbb{C}$ such that the operator $A_k : \mathrm{H}_0^2(\Omega) \to \mathrm{H}_0^2(\Omega)$ introduced in (5) is an isomorphism (we recall that $B = A_0$). The technique we proposed in this paragraph does not allow to obtain this result. We refer the reader to [45] for a proof based on an equivalent formulation of (2). Let us underline that in this paper, the author also requires the assumption that $\sigma$ is uniformly positive or uniformly negative in a neighbourhood of the boundary $\partial \Omega$.

## 3.2 Fredholm property can be lost when $\sigma$ changes sign on the boundary

Now, our goal is to understand what happens if this assumption on $\sigma$ is not satisfied. More precisely, what are the properties of $B$ when $\sigma$ changes sign "on" the boundary $\partial \Omega$? Working with techniques which were developed to study elliptic partial differential equations in non smooth domains, we present occurrences where $B : \mathrm{H}_0^2(\Omega) \to \mathrm{H}_0^2(\Omega)$ is not of Fredholm type.

Let us first introduce the notations. The domain $\Omega \subset \mathbb{R}^2$, is partitioned into two subdomains $\Omega_1$, $\Omega_2$ such that $\overline{\Omega} = \overline{\Omega_1} \cup \overline{\Omega_2}$, $\Omega_1 \cap \Omega_2 = \emptyset$. We assume that $\sigma = \sigma_1$ in $\Omega_1$ and $\sigma = \sigma_2$ in $\Omega_2$, where $\sigma_1 > 0$ and $\sigma_2 < 0$ are two constants. The interface $\Sigma := \overline{\Omega_1} \cap \overline{\Omega_2}$ meets the boundary at exactly two point $O$, $O'$ like in Figure 5. At $O$, $\partial\Omega$ and $\Sigma$ are locally straight lines. Therefore, at this point, the domain $\Omega$, $\Omega_1$ and $\Omega_2$ coincide locally with the unbounded sectors

$$\begin{aligned}
\Xi &:= \{(r\cos\theta, r\sin\theta) \,|\, 0 < r < \infty; \theta \in (0;\pi)\} \\
\Xi_1 &:= \{(r\cos\theta, r\sin\theta) \,|\, 0 < r < \infty; \theta \in (0;\alpha)\} \\
\Xi_2 &:= \{(r\cos\theta, r\sin\theta) \,|\, 0 < r < \infty; \theta \in (\alpha;\pi)\},
\end{aligned}$$

for some $\alpha \in (0;\pi)$.

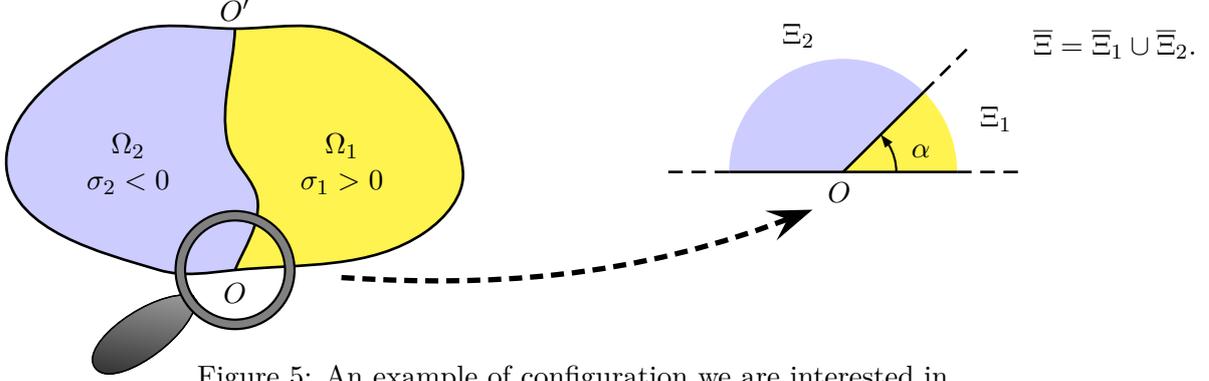

Figure 5: An example of configuration we are interested in.

When one is interested in studying the regularity of the solutions of problem $(\mathscr{P})$, according to [29], we know that a correct start is to compute the *singularities*, that is the non trivial functions of the form

$$s(\boldsymbol{x}) = r^\lambda \varphi(\theta) \tag{34}$$

which satisfy the problem

$$\Delta(\sigma \Delta s) = 0 \text{ a.e. in } \Xi \quad \text{and} \quad s = \partial_\nu s = 0 \text{ a.e. on } \partial \Xi.$$

If $s$ satisfies the above equations, then $\lambda \in \mathbb{C}$ is called a *singular exponent*. We denote $\Lambda_{\alpha,\kappa_\sigma}$ the set of singular exponents.

### 3.2.1 Computation of the singularities

Actually, we will focus our attention only on the computation of some particular singularities: the ones for which the singular exponent verifies $\lambda = 1 + i\eta$, for some $\eta \in \mathbb{R}^* := \mathbb{R} \setminus \{0\}$. The reason is that if such a singularity exists then, as we will prove in §3.2.2, the operator $B : \mathrm{H}_0^2(\Omega) \to \mathrm{H}_0^2(\Omega)$ is not of Fredholm type. We denote $s_1 = s|_{\Xi_1}$, $s_2 = s|_{\Xi_2}$, $\varphi_1 := \varphi|_{(0;\alpha)}$ and $\varphi_2 := \varphi|_{(\alpha;\pi)}$. If $s$ is a singularity, then $(s_1, s_2)$ satisfies the following transmission problem

$$\left|\begin{aligned}
\Delta\Delta s_1 &= 0 & &\text{in } \Xi_1 \\
\Delta\Delta s_2 &= 0 & &\text{in } \Xi_2 \\
s_1 = \nu \cdot \nabla s_1 &= 0 & &\text{on } \partial\Xi \cap \partial\Xi_1 \setminus \{0\} \\
s_2 = \nu \cdot \nabla s_2 &= 0 & &\text{on } \partial\Xi \cap \partial\Xi_2 \setminus \{0\} \\
s_1 - s_2 = \nu_\Sigma \cdot \nabla(s_1 - s_2) &= 0 & &\text{on } \partial\Xi_1 \cap \partial\Xi_2 \setminus \{0\} \\
\sigma_1 \Delta s_1 - \sigma_2 \Delta s_2 = \nu_\Sigma \cdot \nabla(\sigma_1 \Delta s_1 - \sigma_2 \Delta s_2) &= 0 & &\text{on } \partial\Xi_1 \cap \partial\Xi_2 \setminus \{0\}.
\end{aligned}\right.$$

In these equations, $\nu$ (resp. $\nu_\Sigma$) denotes the unit outward normal vector to $\partial\Xi$ (resp. $\Sigma$) oriented to the exterior of $\Xi$ (resp. $\Xi_1$). In polar coordinates, the bilaplacian operator takes the form

$$\Delta^2 = r^{-4}(\partial_\theta^2 + (r\partial_r - 2)^2)(\partial_\theta^2 + (r\partial_r)^2).$$

Imposing the boundary conditions on $\partial\Xi \cap \partial\Xi_1 \setminus \{0\}$ and $\partial\Xi \cap \partial\Xi_2 \setminus \{0\}$, we prove (see [31, §7.1.2]) that for $\lambda \in \mathbb{C} \setminus \{0, 1, 2\}$, the functions $\varphi_1$ and $\varphi_2$ admit the following expressions

$$\varphi_1(\theta) = A\left(\cos(\lambda\theta) - \cos((\lambda-2)\theta)\right) + B\left((\lambda-2)\sin(\lambda\theta) - \lambda\sin((\lambda-2)\theta)\right),$$
$$\begin{aligned}\varphi_2(\theta) = {}&C\left(\cos(\lambda(\theta-\pi)) - \cos((\lambda-2)(\theta-\pi))\right) \\ &+ D\left((\lambda-2)\sin(\lambda(\theta-\pi)) - \lambda\sin((\lambda-2)(\theta-\pi))\right),\end{aligned}$$

where $A$, $B$, $C$ and $D$ are some constants. Writing the transmission conditions at $\theta = \alpha$, we obtain a system of four equations with four unknowns. Computing the determinant, we find that $\lambda = 1 + i\eta$, with $\eta \in \mathbb{R}^*$ is a singular exponent if and only there holds

$$\begin{aligned} 0 = & -1 + \kappa_\sigma - \kappa_\sigma^2 + \eta^2(1-\kappa_\sigma)^2(\cos(2\alpha)-1) \\ & + \kappa_\sigma \cosh(2\pi\eta) + \kappa_\sigma(\kappa_\sigma-1)\cosh(2\alpha\eta) - (\kappa_\sigma-1)\cosh(2\eta(\pi-\alpha)). \end{aligned} \quad (35)$$

In the above equation, the angle $\alpha$ and the contrast $\kappa_\sigma = \sigma_2/\sigma_1$ are two parameters. The question we wish to solve is the following: for a given problem, *i.e.* for given $\alpha \in (0;\pi)$ and $\kappa_\sigma \in (-\infty;0)$, can we find $\eta \in \mathbb{R}^*$ such that (35) is verified. For $\alpha \in (0;\pi)$ and $\kappa_\sigma \in (-\infty;0)$, we define the function $h_{\alpha,\kappa_\sigma}$ such that for all $\eta \in \mathbb{R}$

$$\begin{aligned} h_{\alpha,\kappa_\sigma}(\eta) = & -1 + \kappa_\sigma - \kappa_\sigma^2 + \eta^2(1-\kappa_\sigma)^2(\cos(2\alpha)-1) \\ & + \kappa_\sigma \cosh(2\pi\eta) + \kappa_\sigma(\kappa_\sigma-1)\cosh(2\alpha\eta) - (\kappa_\sigma-1)\cosh(2\eta(\pi-\alpha)). \end{aligned} \quad (36)$$

First, we notice that $h_{\alpha,\kappa_\sigma}$ is even: if $1 + i\eta$, with $\eta \in \mathbb{R}^*$ is a singular exponent, then $1 - i\eta$ is also a singular exponent. Therefore, is it sufficient to study $\eta \mapsto h_{\alpha,\kappa_\sigma}(\eta)$ on $(0;+\infty)$. Then, we observe that there holds $h_{\alpha,\kappa_\sigma}(\eta) = 0$ if and only there holds $h_{\pi-\alpha,1/\kappa_\sigma}(\eta) = 0$. This is reassuring since the singularities for the problem with an angle of aperture $\pi - \alpha$ and a contrast equal to $1/\kappa_\sigma$ are the same as the singularities for the problem with an angle of aperture $\alpha$ and a contrast equal to $\kappa_\sigma$. For all $(\alpha, \kappa_\sigma) \in (0;\pi) \times (-\infty;0)$, we have $h_{\alpha,\kappa_\sigma}(\eta) \to -\infty$ when $\eta \to +\infty$. Moreover, there holds $h_{\alpha,\kappa_\sigma}(0) = 0$. A Taylor expansion of $h_{\alpha,\kappa_\sigma}$ at $\eta = 0$ gives

$$h_{\alpha,\kappa_\sigma}(\eta) = g_\alpha(\kappa_\sigma) \eta^2 + O(\eta^4) \quad (37)$$

with $g_\alpha(\kappa_\sigma) = 2(\alpha^2 - \sin^2(\alpha))\kappa_\sigma^2 - 4(\alpha^2 - \sin^2(\alpha) - \alpha\pi)\kappa_\sigma + 2(\alpha^2 - \sin^2(\alpha) + \pi^2 - 2\alpha\pi)$. For a given $\alpha \in (0;\pi)$, we find that $g_\alpha(\kappa_\sigma)$ is strictly positive when

$$\kappa_\sigma \in I(\alpha) := (-\infty; \ell_-(\alpha)) \cup (\ell_+(\alpha); 0) \quad (38)$$

where $\qquad \ell_-(\alpha) := -\dfrac{\pi - \alpha + \sin(\pi-\alpha)}{\alpha - \sin\alpha} \quad$ and $\quad \ell_+(\alpha) := -\dfrac{\pi - \alpha - \sin(\pi-\alpha)}{\alpha + \sin\alpha}.$

Let us define the region (see Figure 6)

$$\mathscr{R} := \{(\alpha', \kappa_\sigma') \in (0;\pi) \times (-\infty;0) \,|\, \kappa_\sigma' \in I(\alpha')\}. \quad (39)$$

When $(\alpha, \kappa_\sigma)$ belongs to $\mathscr{R}$, the continuous function $h_{\alpha,\kappa_\sigma}$ satisfies $h_{\alpha,\kappa_\sigma}(0) = 0$, $h_{\alpha,\kappa_\sigma}(\eta) = g_\alpha(\kappa_\sigma) \eta^2 + O(\eta^4)$ at $\eta = 0$, with $g_\alpha(\kappa_\sigma) > 0$, and $\lim_{\eta \to +\infty} h_{\alpha,\kappa_\sigma}(\eta) = -\infty$. This allows to deduce that $h_{\alpha,\kappa_\sigma}$ vanishes at least once on $(0;+\infty)$. Thus, if $(\alpha, \kappa_\sigma) \in \mathscr{R}$, then there exist singularities of the form $s(r,\theta) = r^{1+i\eta}\varphi(\theta)$ with $\eta \in \mathbb{R}^*$. We will complement these results in Section 5. In Proposition 5.1, we will prove that if $(\alpha, \kappa_\sigma) \in \mathscr{R}$, then the set of singular exponents $\Lambda_{\alpha,\kappa_\sigma}$ is such that $\Lambda_{\alpha,\kappa_\sigma} \cap \{\lambda \in \mathbb{C} \setminus \{1\} \,|\, \Re e\, \lambda = 1\} = \{1 \pm i\eta_0\}$, for some $\eta_0 > 0$. On the other hand, in Proposition 5.2, we will establish that if $(\alpha, \kappa_\sigma)$ is located in $((0;\pi) \times (-\infty;0)) \setminus \overline{\mathscr{R}}$, then $\Lambda_{\alpha,\kappa_\sigma} \cap \{\lambda \in \mathbb{C} \setminus \{1\} \,|\, \Re e\, \lambda = 1\} = \emptyset$.

### 3.2.2 Proof of ill-posedness

Now, we can establish the

**Proposition 3.1.** *Assume that $(\sigma, \kappa_\sigma)$ belongs to the region $\mathscr{R}$ defined in (39). Then, the operator $B : \mathrm{H}_0^2(\Omega) \to \mathrm{H}_0^2(\Omega)$ defined in (7) is not of Fredholm type.*

*Proof.* Our proof relies on the Peetre's lemma [42] (see also lemma 5.1 in [33, Chap. 2], or lemma 3.4.1 in [30]):

**Lemma 3.2.** *Let* X*,* Y *and* Z *be three Banach spaces such that* X *is compactly embedded into* Z*. Let* $L : X \to Y$ *be a continuous linear map. Then the assertions below are equivalent:*

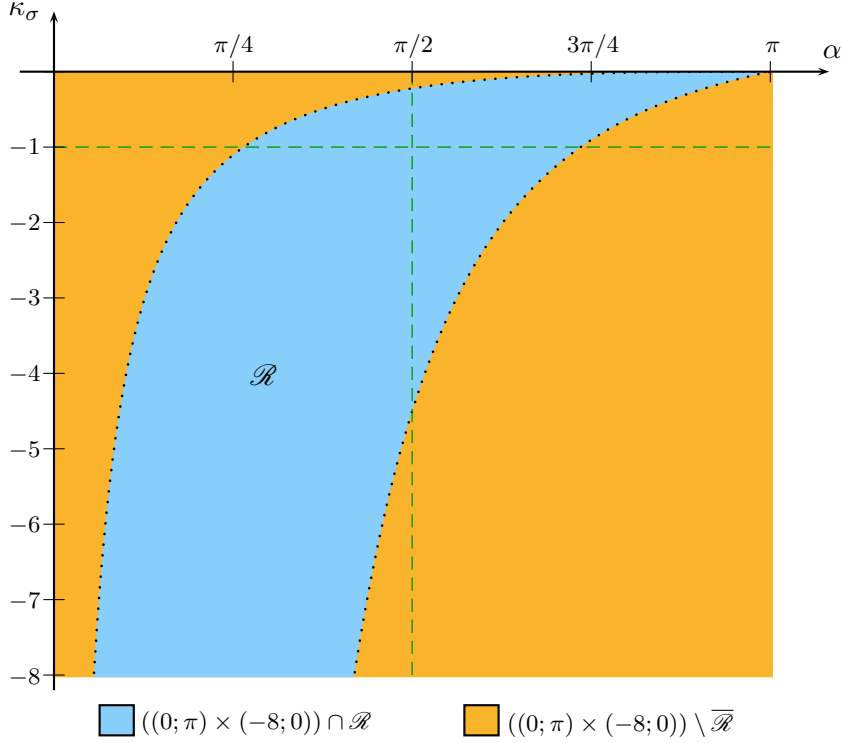

Figure 6: When $(\alpha, \kappa_\sigma)$ belongs to the orange region $\mathscr{R}$, there exist singularities of the form $s(r,\theta) = r^{1+i\eta}\varphi(\theta)$ with $\eta \in \mathbb{R}^*$. Therefore, in these situations the operator $B : \mathrm{H}_0^2(\Omega) \to \mathrm{H}_0^2(\Omega)$ defined in (7) is not of Fredholm type.

> i) $\dim(\ker L) < \infty$ *and* range $L$ *is closed in* $\mathrm{Y}$;
> ii) *there exists* $C > 0$ *such that*
> $$\|v\|_\mathrm{X} \leq C(\,\|Lv\|_\mathrm{Y} + \|v\|_\mathrm{Z}\,), \qquad \forall v \in \mathrm{X}.$$

When $(\sigma, \kappa_\sigma)$ belongs to the region $\mathscr{R}$, we know that there exists a singularity $s(\boldsymbol{x}) = r^{1+i\eta}\varphi(\theta)$ with $\eta \in \mathbb{R}^*$ such that $\Delta(\sigma\Delta s) = 0$ a.e. in $\Xi$ and $s = \partial_\nu s = 0$ a.e. on $\partial \Xi$. Let $\zeta \in \mathscr{C}_0^\infty([0;+\infty), [0;1])$ be a cut-off function equal to 1 in a neighbourhood of $O$. For $m \in \mathbb{N}^* = \{1,2,3,\dots\}$, we define

$$s_m(\boldsymbol{x}) := r^{1+i\eta+1/m}\varphi(\theta) \qquad \text{and} \qquad v_m(\boldsymbol{x}) := \zeta(r)s_m(\boldsymbol{x}).$$

We assume that the support of $\zeta$ is such that $v_m|_{\partial \Omega} = 0$. Since $\Re\mathrm{e}\,(1+i\eta+1/m) > 1$, one can check by a direct computation that the function $v_m$ belongs to $\mathrm{H}_0^2(\Omega)$ for all $m \in \mathbb{N}^*$. Our goal is to establish the following properties

$$\lim_{m\to+\infty} \|v_m\|_{\mathrm{H}_0^2(\Omega)} = +\infty \qquad \text{and} \qquad \|Bv_m\|_{\mathrm{H}^{-2}(\Omega)} + \|v_m\|_{\mathrm{H}_0^1(\Omega)} \leq C, \qquad \forall m \in \mathbb{N}^*.$$

Together with Lemma 3.2, this will prove that $B$ is not of Fredholm type when $(\sigma, \kappa_\sigma) \in \mathscr{R}$.

★ Behaviour of $(\|v_m\|_{\mathrm{H}_0^2(\Omega)})_m$. By definition, we have $\|v_m\|_{\mathrm{H}_0^2(\Omega)} = \|\Delta v_m\|_\Omega \geq \|\Delta v_m\|_{\widetilde{\Omega}}$, where $\widetilde{\Omega} := \{\boldsymbol{x} \in \Omega \,|\, \zeta(r) = 1\}$. On $\widetilde{\Omega}$, we find

$$\Delta v_m = r^{-1+i\eta+1/m}((1+i\eta+1/m)^2 \varphi(\theta) + \varphi''(\theta)).$$

Let us define the function $\varsigma_m$ such that $\varsigma_m(\theta) = (1+i\eta+1/m)^2\varphi(\theta) + \varphi''(\theta)$. We have

$$\begin{aligned}\|\Delta v_m\|_{\widetilde{\Omega}}^2 &\geq \int_0^\delta \int_0^\pi r^{-2+2/m}|\varsigma_m(\theta)|^2 \, r d\theta dr \\ &\geq \|\varsigma_m\|_{(0;\pi)}^2 \int_0^\delta r^{-1+2/m} dr = \|\varsigma_m\|_{(0;\pi)}^2 \, \frac{m}{2}\, \delta^{2/m}.\end{aligned} \qquad (40)$$

Above $\delta > 0$ is a fixed small number such that $\zeta(r) = 1$ on $[0;\delta]$. Since $(1+i\eta)^2\varphi + \varphi'' \not\equiv 0$ (otherwise, we should have $\Delta s = 0$, which is impossible due to the boundary conditions), there holds $\|\varsigma_m\|^2_{(0;\pi)} \neq 0$ for $m$ large enough. From (40), we deduce that $\|v_m\|_{H^2_0(\Omega)} \underset{m \to +\infty}{\to} +\infty$.

★ Behaviour of $(\|v_m\|_{H^1_0(\Omega)})_m$. By a direct computation, we can check that there exists a constant $C > 0$ such that for all $m \in \mathbb{N}^*$, we have $\|v_m\|_{H^1_0(\Omega)} \leq C$.

★ Behaviour of $(\|Bv_m\|_{H^{-2}(\Omega)})_m$. Now, let us prove that the sequence $(\|\Delta(\sigma\Delta v_m)\|_{H^{-2}(\Omega)})_m$ remains bounded. By definition, we have

$$\|\Delta(\sigma\Delta v_m)\|_{H^{-2}(\Omega)} := \sup_{v \in H^2_0(\Omega)\,|\,\|\Delta v\|_\Omega = 1} |(\sigma\Delta v_m, \Delta v)_\Omega| = \sup_{v \in \mathscr{C}_0^\infty(\Omega)\,|\,\|\Delta v\|_\Omega = 1} |(\sigma\Delta v_m, \Delta v)_\Omega|.$$

We compute, for all $v \in \mathscr{C}_0^\infty(\Omega)$, $\Delta v_m = \Delta(\zeta s_m) = \zeta \Delta s_m + 2\nabla s_m \nabla \zeta + s_m \Delta \zeta$ and $\Delta(\zeta v) = \zeta \Delta v + 2\nabla v \nabla \zeta + v \Delta \zeta$. This allows to write

$$\begin{aligned}(\sigma\Delta v_m, \Delta v)_\Omega &= (\sigma(\zeta \Delta s_m + 2\nabla v_m \nabla \zeta + v_m \Delta \zeta), \Delta v)_\Omega \\ &= (\sigma\Delta s_m, \Delta(\zeta v) - 2\nabla v \nabla \zeta - v\Delta \zeta)_\Omega + (\sigma(2\nabla s_m \nabla \zeta + s_m \Delta \zeta), \Delta v)_\Omega.\end{aligned}$$

We deduce

$$\begin{aligned}|(\sigma\Delta v_m, \Delta v)_\Omega| &\leq \underbrace{|(\sigma\Delta s_m, \Delta(\zeta v))_\Omega|}_{\text{❶}} + 2\underbrace{|(\sigma\Delta s_m, \nabla v \nabla \zeta)_\Omega|}_{\text{❷}} + \underbrace{|(\sigma\Delta s_m, v\Delta\zeta)_\Omega|}_{\text{❸}} \\ &\quad + 2\underbrace{|(\sigma\nabla s_m \nabla \zeta, \Delta v)_\Omega|}_{\text{❹}} + \underbrace{|(\sigma s_m \Delta\zeta, \Delta v)_\Omega|}_{\text{❺}}.\end{aligned}$$

Let us study the term ❶. Integrating by parts, we find $(\sigma\Delta s_m, \Delta(\zeta v))_\Omega = (\Delta(\sigma\Delta s_m), \zeta v)_\Omega$. A direct computation allows to establish that $\Delta(\sigma\Delta s_m)(\boldsymbol{x}) = r^{-3+i\eta+1/m}\tilde{\varphi}_m(\theta)$ where $\tilde{\varphi}_m \in L^2(0;\pi)$ is such that $\|\tilde{\varphi}_m\|_{(0;\pi)} \leq C/m$ for some constant $C$ independent of $m$. Therefore, we can write

$$\begin{aligned}|(\sigma\Delta s_m, \Delta(\zeta v))_\Omega| &= |(r^{-3+i\eta+1/m}\tilde{\varphi}_m, \zeta v)_\Omega| \\ &= C\,|(r^{-1+i\eta+1/m}\tilde{\varphi}_m, r^{-2}(r\partial_r)^2(\zeta v))_\Omega| \\ &\leq C\,\|r^{-1+i\eta+1/m}\tilde{\varphi}_m\|_\Omega \|v\|_{H^2_0(\Omega)} \leq C^{2/m}\|v\|_{H^2_0(\Omega)}.\end{aligned}$$

In the above estimate, $C$ is a constant independent of $m$ which can change from one line to another. The second line is obtained proceeding to an integration by part with respect to the $r$ variable. This proves that ❶ is bounded as $m \to +\infty$. To deal with the terms ❷ and ❸, we use that $\nabla \zeta$ and $\Delta \zeta$ vanish in a neighbourhood of $O$ to obtain

$$|(\sigma\Delta s_m, \nabla v \nabla \zeta)_\Omega| = |(\sigma\Delta s_m, \nabla v \nabla \zeta)_{\Omega \setminus B(O,\delta)}| \leq C\,\|\Delta s_m\|_{\Omega \setminus B(O,\delta)} \|v\|_{H^2_0(\Omega)}$$
$$\text{and}\quad |(\sigma\Delta s_m, v\Delta\zeta)_\Omega| = |(\sigma\Delta s_m, v\Delta\zeta)_{\Omega\setminus B(O,\delta)}| \leq C\,\|\Delta s_m\|_{\Omega\setminus B(O,\delta)} \|v\|_{H^2_0(\Omega)}.$$

Since the sequence $(\|\Delta v_m\|_{\Omega \setminus B(O,\delta)})_m$ is bounded, we deduce that the terms ❷ and ❸ remain bounded. Finally, to deal with the terms ❹ and ❺, we write

$$|(\sigma\nabla s_m \nabla \zeta, \Delta v)_\Omega| + |(\sigma s_m \Delta \zeta, \Delta v)_\Omega| \leq C\,\|s_m\|_{H^1_0(\Omega)} \|v\|_{H^2_0(\Omega)}.$$

Since the sequence $(\|s_m\|_{H^1_0(\Omega)})_m$ is bounded, we deduce that the terms ❹ and ❺ are bounded. All these intermediate results allow to conclude that the sequence $(\|\Delta(\sigma\Delta v_m)\|_{H^{-2}(\Omega)})_m$ remains bounded as $m \to +\infty$. $\square$

# 4 Appendix: bilaplacian with mixed boundary conditions in dimension $d \geq 3$

In Section 2, we studied the properties of problem $(\tilde{\mathscr{P}})$, defined in (10), for which we impose mixed boundary conditions. In Theorem 2.1, we proved that $(\tilde{\mathscr{P}})$ is well-posed when $\sigma$ satisfies $\sigma \in L^\infty(\Omega)$, $\sigma^{-1} \in L^\infty(\Omega)$ and when the domain $\Omega$ is smooth or convex. For 2D polygons with reentrant corners, we established in Propositions 2.5 and 2.9 that a kernel and a cokernel can appear for the operator $\tilde{B}$ associated with $(\tilde{\mathscr{P}})$. In this section, we present some results for $(\tilde{\mathscr{P}})$ when $\Omega \subset \mathbb{R}^d$, $d \geq 3$, is neither smooth nor convex. As we saw previously, the features of $(\tilde{\mathscr{P}})$ are governed by the properties of the operator $\Delta : H_0^1(\Omega) \cap H^2(\Omega) \to L^2(\Omega)$. In 2D, for an open set with a polygonal boundary, we have the following results. Either the domain $\Omega$ is convex and then $\Delta : H_0^1(\Omega) \cap H^2(\Omega) \to L^2(\Omega)$ is an isomorphism. Or, $\Omega$ has some reentrant corners and then $\Delta : H_0^1(\Omega) \cap H^2(\Omega) \to L^2(\Omega)$ is an injective Fredholm operator with a cokernel whose dimension is equal to the number of reentrant corners. In arbitrary dimension $d \geq 3$, the features of $\Delta : H_0^1(\Omega) \cap H^2(\Omega) \to L^2(\Omega)$ are more varied.

### 4.1 Non convex conical tips in dimension $d$

First, let us focus our attention on domains whose boundary has concave *conical tips*. To simplify the presentation, we assume that $\partial\Omega$ possesses only one such singularity: the boundary of $\Omega \subset \mathbb{R}^d$ is of class $\mathscr{C}^2$ except at $O$. At this point, $\partial\Omega$ coincides locally with a cone. Let us clarify this notion. Consider $\omega$ a domain of $\mathbb{S}^{d-1}$, the unit sphere of $\mathbb{R}^d$. Define $K_\omega^R := \{r\boldsymbol{\theta} \,|\, 0 < r < R,\, \boldsymbol{\theta} \in \omega\}$. We assume that there exists $R > 0$ and $\omega \subset \mathbb{S}^{d-1}$ of class $\mathscr{C}^2$ such that $\Omega \cap B(R, O) = K_\omega^R$. Here, we define $B(R, O) := \{\boldsymbol{x} \in \mathbb{R}^d \,|\, |\boldsymbol{x}| < R\}$.

By virtue of the Lax-Milgram theorem, the operator $\Delta : H_0^1(\Omega) \cap H^2(\Omega) \to L^2(\Omega)$ is injective whatever the dimension of space $d$. Hence, the main task consists in checking whether or not its range is closed and in studying its cokernel if the answer to the latter question is positive. To address such a problem, according to [29], we know that it is reasonable to start by computing the non trivial functions of the form

$$s(\boldsymbol{x}) = r^\Lambda \Phi(\boldsymbol{\theta}) \tag{41}$$

which satisfy the problem

$$\Delta s = 0 \text{ a.e. in } K_\omega^\infty \qquad \text{and} \qquad s = 0 \text{ a.e. on } \partial K_\omega^\infty. \tag{42}$$

Here, we denote $K_\omega^\infty := \{r\boldsymbol{\theta} \,|\, 0 < r,\, \boldsymbol{\theta} \in \omega\}$. In polar coordinates, the Laplace operator is given by the formula

$$\Delta = \frac{\partial^2}{\partial r^2} + \frac{d-1}{r}\frac{\partial}{\partial r} + \frac{1}{r^2}\tilde{\Delta},$$

where $\tilde{\Delta}$ denotes the Laplace-Beltrami operator on the unit sphere. The singular exponents in (41), *i.e.* the values of $\Lambda$ for which $s$ satisfies Equations (42), are defined by:

$$\Lambda_k^\pm := 1 - \frac{d}{2} \pm \sqrt{\left(1 - \frac{d}{2}\right)^2 + \mu_k}. \tag{43}$$

Above, $\mu_k$ is the $k$th eigenvalue of the problem

$$\left|\begin{array}{rcll} -\tilde{\Delta}\Phi(\boldsymbol{\theta}) &=& \mu\,\Phi(\boldsymbol{\theta}) & \text{in } \omega \\ \Phi(\boldsymbol{\theta}) &=& 0 & \text{on } \partial\omega \end{array}\right., \tag{44}$$

and the functions $\Phi$ in (41) are equal to the eigenvectors of problem (44). The Laplace-Beltrami operator is selfadjoint and positive definite. Therefore, the eigenvalues of Problem (44) form the sequence

$$0 < \mu_1 < \mu_2 \leq \mu_3 \leq \ldots \qquad \text{with } \mu_k \to +\infty \text{ when } k \to +\infty.$$

Classically, the first eigenvalue $\mu_1$ is simple (cf. [25, theorem 1.2.5]). Therefore, the positive singular exponents in (41) satisfy

$$0 < \Lambda_1^+ < \Lambda_2^+ \leq \Lambda_3^+ \leq \ldots \qquad \text{with } \Lambda_k^+ \to +\infty \text{ when } k \to +\infty.$$

The negative singular exponents are given by the expression $\Lambda_k^- = 2 - d - \Lambda_k^+$.

Now, we introduce well-suited spaces to handle the singularities (41). For $l \in \mathbb{N}$ and $\beta \in \mathbb{R}$, we define the space $\mathrm{V}_\beta^l(\Omega)$ as the closure of $\mathscr{C}_0^\infty(\overline{\Omega}\setminus\{O\}) := \{\varphi \in \mathscr{C}^\infty(\overline{\Omega}) \,|\, [\exists \delta > 0 \,|\, \mathrm{supp}(\varphi) \cap B(O, \delta) = \emptyset]\}$ for the norm

$$\|\varphi\|_{\mathrm{V}_\beta^l(\Omega)} := \Big(\sum_{|\boldsymbol{\alpha}| \leq l} \int_\Omega r^{2(\beta - l + |\boldsymbol{\alpha}|)} |\partial_{\boldsymbol{x}}^{\boldsymbol{\alpha}} \varphi|^2 \, d\boldsymbol{x}\Big)^{1/2}. \tag{45}$$

To take into account the homogeneous Dirichlet boundary condition, we introduce, for $l \in \mathbb{N}^* = \{1, 2, 3, \dots\}$ and $\beta \in \mathbb{R}$, the space $\mathring{\mathrm{V}}_\beta^l(\Omega)$, which is the closure of $\mathscr{C}_0^\infty(\Omega)$ for the norm (45). For $l \in \mathbb{N}^*$ and $\beta \in \mathbb{R}$, we define the continuous operators

$$\begin{aligned} A_\beta^l : \ \mathrm{V}_\beta^{l+1}(\Omega) \cap \mathring{\mathrm{V}}_{\beta - l}^1(\Omega) &\to \mathrm{V}_\beta^{l-1}(\Omega) \\ \varphi &\mapsto A_\beta^l \varphi = \Delta\varphi. \end{aligned} \tag{46}$$

There holds the fundamental theorem (see for example [29, 30, 36, 35, 17, 18]):

**Theorem 4.1.** *Let $\Lambda_1^+$ be defined in (43). For $l \in \mathbb{N}^*$ and $\beta \in \mathbb{R}$, the operator $A_\beta^l : \mathrm{V}_\beta^{l+1}(\Omega) \cap \mathring{\mathrm{V}}_{\beta - l}^1(\Omega) \to \mathrm{V}_\beta^{l-1}(\Omega)$ is an isomorphism if and only if*

$$1 - \Lambda_1^+ < \beta - l + \frac{d}{2} < d - 1 + \Lambda_1^+. \tag{47}$$

*More precisely,*
*i) if $\beta - l + d/2 < 1 - \Lambda_1^+$, $A_\beta^l$ is of Fredholm type injective but not onto;*
*ii) if $\beta - l + d/2 > d - 1 + \Lambda_1^+$, $A_\beta^l$ is of Fredholm type onto but not injective;*
*iii) if $\beta - l + d/2 = 1 - \Lambda_1^+$ or if $\beta - l + d/2 = d - 1 + \Lambda_1^+$, the operator $A_\beta^l$ is not of Fredholm type because its range is not closed in $\mathrm{V}_\beta^{l-1}(\Omega)$.*

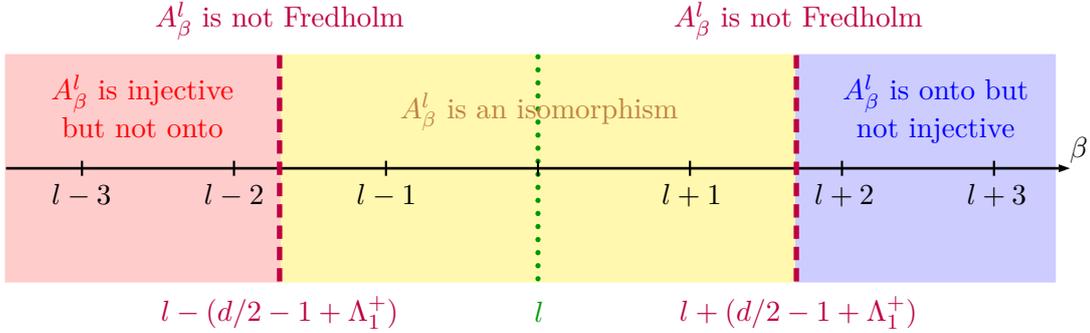

Figure 7: Features of the operator $A_\beta^l : \mathrm{V}_\beta^{l+1}(\Omega) \cap \mathring{\mathrm{V}}_{\beta - l}^1(\Omega) \to \mathrm{V}_\beta^{l-1}(\Omega)$ with respect to $\beta$, $d$ being the dimension of the space.

For the study of our problem, the properties of the operator $A_0^1$, i.e. $A_\beta^l$ with $l = 1$ and $\beta = 0$, will be useful. Indeed, we notice that $\mathrm{V}_0^0(\Omega) = \mathrm{L}^2(\Omega)$. Moreover, [38, lemma 3.4] establishes that $\mathrm{V}_0^2(\Omega) \cap \mathring{\mathrm{V}}_{-1}^1(\Omega) = \mathrm{H}_0^1(\Omega) \cap \mathrm{H}^2(\Omega)$.

<div align="center">

★ ★ ★ ★ ★

Conical tips in dimension $d \geq 4$

★ ★ ★ ★ ★

</div>

Theorem 4.1 indicates that the operator $A_0^1 : V_0^2(\Omega) \cap \overset{\circ}{V}{}_{-1}^1(\Omega) \to V_0^0(\Omega)$ is an isomorphism if and only if
$$1 - \Lambda_1^+ < 0 - 1 + d/2 < d - 1 + \Lambda_1^+ \quad \Leftrightarrow \quad d > 4 - 2\Lambda_1^+.$$

This is always true in dimension $d \geq 4$. Since $V_0^2(\Omega) \cap \overset{\circ}{V}{}_{-1}^1(\Omega) = H_0^1(\Omega) \cap H^2(\Omega)$, we deduce that the operator $\Delta : H_0^1(\Omega) \cap H^2(\Omega) \to L^2(\Omega)$ is an isomorphism. Proceeding as in the proof of Theorem 2.1, with the T-coercivity technique, we show easily the

**Proposition 4.1.** *Let $\Omega \subset \mathbb{R}^d$, $d \geq 4$, be a domain whose boundary is of class $\mathscr{C}^2$ except at a finite number of points where it coincides locally with a non convex conical tip. Assume that $\sigma \in L^\infty(\Omega)$ is such that $\sigma^{-1} \in L^\infty(\Omega)$. Then the operator $\tilde{B} : H_0^1(\Omega) \cap H^2(\Omega) \to H_0^1(\Omega) \cap H^2(\Omega)$ defined in (12) is an isomorphism.*

$\star\star\star\star\star$

Conical tips in dimension $d = 3$

$\star\star\star\star\star$

In dimension $d = 3$, the relation $d > 4 - 2\Lambda_1^+$ is not always true. This depends on the value of $\Lambda_1^+$, that is on the value of $\mu_1$, the first eigenvalue of the Laplace-Beltrami operator on the portion of sphere $\omega$.

The eigenvalue $\mu_1(\omega)$ depends continuously on the domain $\omega$ (see [27, 25]) and, according to the *min-max* principle, $\omega_a \subset \omega_b$ implies $\mu_1(\omega_a) \geq \mu_1(\omega_b)$. Therefore, when $\Omega$ is convex, we have $\mu_1(\omega) > \mu_1((\mathbb{R}^2 \times \mathbb{R}_+^*) \cap \mathbb{S}^2) = 2$ and so $\Lambda_1^+ > 1$ (here, $\mathbb{R}_+^* = (0; +\infty)$). Thus, when the conical tip is convex, the relation $d > 4 - 2\Lambda_1^+$ is verified as soon as $d \geq 2$ and the operator $\Delta : H_0^1(\Omega) \cap H^2(\Omega) \to L^2(\Omega)$ is an isomorphism. This is coherent with the theorem 3.2.1.2 of [23] we have used which indicates that $\Delta : H_0^1(\Omega) \cap H^2(\Omega) \to L^2(\Omega)$ is an isomorphism in any dimension as soon as the domain $\Omega$ is convex. This result remains valid when the conical tip is non convex with $\Lambda_1^+ > 1/2 \Leftrightarrow \mu_1 > 3/4$. This allows us to state the

**Proposition 4.2.** *Let $\Omega \subset \mathbb{R}^3$ be a domain whose boundary is of class $\mathscr{C}^2$ except at one point $O$ where it coincides with a conical tip of aperture $\omega \subset \mathbb{S}^2$. Assume that $\omega$ is such that $\mu_1(\omega)$, the first eigenvalue of the Laplace-Beltrami operator defined in (44), verifies $\mu_1(\omega) > 3/4$. Assume that $\sigma \in L^\infty(\Omega)$ is such that $\sigma^{-1} \in L^\infty(\Omega)$. Then, the operator $\tilde{B} : H_0^1(\Omega) \cap H^2(\Omega) \to H_0^1(\Omega) \cap H^2(\Omega)$ defined in (12) is an isomorphism.*

**Remark 4.1.** *Let us emphasize that the result of this proposition is not empty! Indeed, there exist non convex conical tips for which $\mu_1(\omega) > 3/4$. To be convinced of this, it is sufficient to remember that $\mu_1((\mathbb{R}^2 \times \mathbb{R}_+^*) \cap \mathbb{S}^2) = 2$ and that $\mu_1(\omega)$ depends continuously on $\omega$.*

In [38, lemma 5.1], the authors exhibit conical tips such that $\mu_1(\omega) < 3/4$. This allows them to deduce that there exist conical tips such that $\mu_1(\omega) = 3/4$ (see also [16, 12, 1]). For the latter ones, according to Theorem 4.1, the operator $\Delta : H_0^1(\Omega) \cap H^2(\Omega) \to L^2(\Omega)$ is not of Fredholm type because its range is not closed in $L^2(\Omega)$. For such configurations, we cannot use the process of resolution in two steps to prove that $(\tilde{\mathscr{P}})$ is well-posed and the example 5.3.2 of [38] show that the operator $\tilde{B}$ is not of Fredholm type for $\sigma = 1$. In the sequel, we shall discard this case and assume that $\Omega$ is such that $\mu_1(\omega) < 3/4$.

The lemma 5.2 of [38] indicates that $\Lambda_2^+ > 1$. As a consequence, when $\mu_1(\omega) \in (0; 3/4)$, the cokernel of $\Delta : H_0^1(\Omega) \cap H^2(\Omega) \to L^2(\Omega)$ is of dimension equal to one and we can use again the technique of the proofs of Propositions 2.8 and 2.9. To do so, let us introduce the function $\zeta$ such that
$$\zeta(\boldsymbol{x}) = r^{-1-\Lambda_1^+} \Phi_1(\boldsymbol{\theta}) + \tilde{\zeta}(\boldsymbol{x}), \tag{48}$$

where $\tilde\zeta$ is the unique element of $\mathrm{H}^1(\Omega)$ verifying $\Delta\tilde\zeta = 0$ a.e. in $\Omega$ and $\tilde\zeta = -r^{-1-\Lambda_1^+}\Phi_1(\boldsymbol{\theta})$ a.e. on $\partial\Omega$. Next, we consider the function $\psi \in \mathrm{H}_0^1(\Omega)$ which satisfies

$$\Delta\psi = \sigma^{-1}\zeta \qquad (49)$$

and we define $\xi \in \mathrm{H}_0^1(\Omega)$ such that

$$(\nabla\xi, \nabla\xi')_\Omega = (\sigma^{-1}\zeta, \xi')_\Omega, \qquad \forall \xi' \in \mathrm{H}_0^1(\Omega). \qquad (50)$$

Now, we can state the

**Proposition 4.3.** *Assume that $\Omega \subset \mathbb{R}^3$ is a domain whose boundary is of class $\mathscr{C}^2$ except at one point $O$ where it coincides with a cone of aperture $\omega \subset \mathbb{S}^2$. Assume that $\omega$ is such that $\mu_1(\omega)$, the first eigenvalue of the Laplace-Beltrami operator defined in ([44](#)), verifies $\mu_1(\omega) < 3/4$.*
• *If $\sigma \in \mathrm{L}^\infty(\Omega)$ satisfies $\sigma^{-1} \in \mathrm{L}^\infty(\Omega)$ and $(\sigma^{-1}\zeta, \zeta)_\Omega \neq 0$, then the operator $\tilde{B} : \mathrm{H}_0^1(\Omega) \cap \mathrm{H}^2(\Omega) \to \mathrm{H}_0^1(\Omega) \cap \mathrm{H}^2(\Omega)$ defined in ([12](#)) is an isomorphism.*
• *If $\sigma \in \mathrm{L}^\infty(\Omega)$ satisfies $\sigma^{-1} \in \mathrm{L}^\infty(\Omega)$ and $(\sigma^{-1}\zeta, \zeta)_\Omega = 0$, then $\ker \tilde{B} = \mathrm{span}(\psi)$ and for all $f \in (\mathrm{H}_0^1(\Omega) \cap \mathrm{H}^2(\Omega))^*$, $(\tilde{\mathscr{P}})$ has a solution if and only if $\langle f, \psi\rangle_\Omega = 0$.*
*In this statement, the functions $\zeta$, $\xi$ and $\psi$ are respectively defined in ([48](#)), ([50](#)) and ([49](#)).*

⋄ EXAMPLE. We introduce the spherical coordinates defined by the relations

$$(x,y,z) = (r\cos\theta, r\sin\theta\cos\varphi, r\sin\theta\sin\varphi).$$

For $R > 0$ and $0 < \alpha < \pi$, we consider the domain

$$\Omega := \{(r\cos\theta, r\sin\theta\cos\varphi, r\sin\theta\sin\varphi),\ 0 < r < R,\ 0 \le \theta < \alpha,\ 0 \le \varphi < 2\pi\}.$$

For $\alpha \in\ ]\pi/2; \pi[$, we have $\zeta(\boldsymbol{x}) = r^{-1-\Lambda_1^+}\Phi_1(\boldsymbol{\theta}) - R^{-1-\Lambda_1^+}(r/R)^{\Lambda_1^+}\Phi_1(\boldsymbol{\theta})$. Let us call $\alpha_c > \pi/2$ the value of $\alpha$ for which there holds $\mu_1(\omega) = 3/4$.

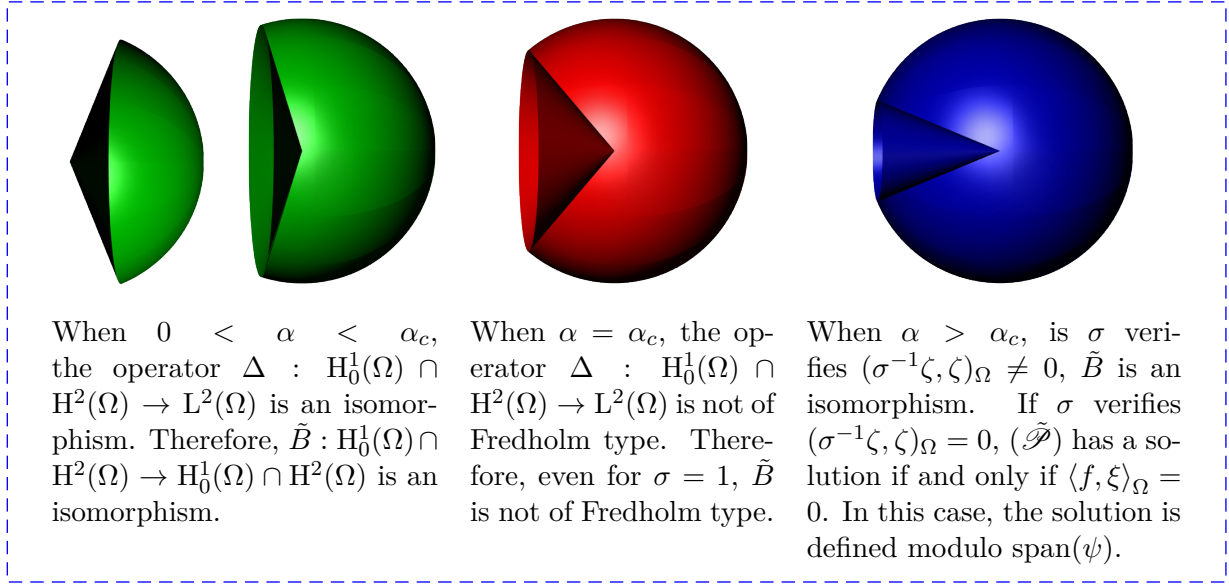

When $0 < \alpha < \alpha_c$, the operator $\Delta : \mathrm{H}_0^1(\Omega) \cap \mathrm{H}^2(\Omega) \to \mathrm{L}^2(\Omega)$ is an isomorphism. Therefore, $\tilde{B} : \mathrm{H}_0^1(\Omega) \cap \mathrm{H}^2(\Omega) \to \mathrm{H}_0^1(\Omega) \cap \mathrm{H}^2(\Omega)$ is an isomorphism.

When $\alpha = \alpha_c$, the operator $\Delta : \mathrm{H}_0^1(\Omega) \cap \mathrm{H}^2(\Omega) \to \mathrm{L}^2(\Omega)$ is not of Fredholm type. Therefore, even for $\sigma = 1$, $\tilde{B}$ is not of Fredholm type.

When $\alpha > \alpha_c$, is $\sigma$ verifies $(\sigma^{-1}\zeta, \zeta)_\Omega \neq 0$, $\tilde{B}$ is an isomorphism. If $\sigma$ verifies $(\sigma^{-1}\zeta, \zeta)_\Omega = 0$, $(\tilde{\mathscr{P}})$ has a solution if and only if $\langle f, \xi\rangle_\Omega = 0$. In this case, the solution is defined modulo $\mathrm{span}(\psi)$.

## 4.2 Non convex edges in dimension $d = 3$

The case of non convex edges in dimension $d = 3$ (see Figure [8](#)) presents additional difficulties. Let us comment on this. For such domains, the injective operator $\Delta : \mathrm{H}_0^1(\Omega) \cap \mathrm{H}^2(\Omega) \to \mathrm{L}^2(\Omega)$ is not of Fredholm type. To make short, the *coefficients* in front of the singularities for conical tips are replaced by *functions* for edges. The cokernel of $\Delta$, equal to $\mathrm{span}(\zeta_1, \ldots, \zeta_N)$ in 2D (cf. Proposition [2.6](#)) is now a space of infinite dimension. When $\sigma = 1$, in [37], the authors achieve to prove that $(\tilde{\mathscr{P}})$ is well-posed extending the method used to consider the 2D case. This times, the compatibility conditions are written for a functional space. We can imagine to implement the same technique to study $(\tilde{\mathscr{P}})$. This would lead to a result of the type: $(\tilde{\mathscr{P}})$ is well-posed if $\sigma$ satisfies an infinite set of "orthogonality" relations. When this condition is not met, we can imagine configurations where the operator $\tilde{B}$ associated with $(\tilde{\mathscr{P}})$ is not of Fredholm type.

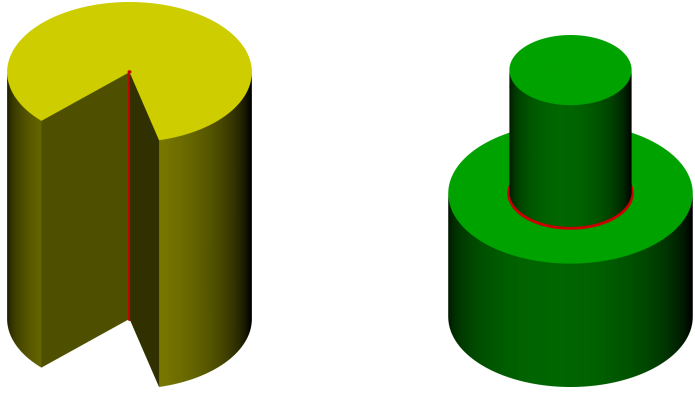

Figure 8: Examples of non convex edges (in red).

## 5 Appendix: complementary computations

In this section, we complement the results we obtained in §3.2.1. We recall that the goal was to compute the singularities of the operator $B$ in a situation where $\sigma$ changes sign on the boundary of $\Omega \subset \mathbb{R}^2$. We use the same notations. Moreover, for $m \in \mathbb{N}$, we call $h^{(m)}_{\alpha,\kappa_\sigma}(\eta_0)$ the $m$th derivative of $\eta \mapsto h_{\alpha,\kappa_\sigma}(\eta)$ at $\eta_0$.

**Proposition 5.1.** *Assume that $(\alpha, \kappa_\sigma) \in \mathscr{R}$ where $\mathscr{R}$ is the region defined in (39). Then, the set of singular exponents $\Lambda_{\alpha,\kappa_\sigma}$ is such that $\Lambda_{\alpha,\kappa_\sigma} \cap \{\lambda \in \mathbb{C} \setminus \{1\} \,|\, \Re e\,\lambda = 1\} = \{1 \pm i\eta_0\}$, for some $\eta_0 > 0$.*

*Proof.* According to the results obtained in §3.2.1, we know that to prove this proposition, it is sufficient to show that the function $\eta \mapsto h_{\alpha,\kappa_\sigma}(\eta)$ defined in (36) vanishes exactly once on $(0; +\infty)$. First, notice that for all $m \in \mathbb{N}$, the continuous function $\eta \mapsto h^{(m)}_{\alpha,\kappa_\sigma}(\eta)$ is such that $h^{(m)}_{\alpha,\kappa_\sigma}(\eta) \underset{\eta \to +\infty}{\sim} \kappa_\sigma(2\pi)^m \cosh(2\pi\eta) \underset{\eta \to +\infty}{\to} -\infty$. This will be a useful property in the analysis.

Let us denote $k_0 \in \{2, 3, 4, \dots\}$ the first number such that $h^{(2(k_0-1))}_{\alpha,\kappa_\sigma}(0) > 0$ and $h^{(2k_0)}_{\alpha,\kappa_\sigma}(0) \leq 0$. This $k_0$ is well-defined since there hold $h^{(2)}_{\alpha,\kappa_\sigma}(0) = g_\alpha(\kappa_\sigma) > 0$ when $(\alpha, \kappa_\sigma) \in \mathscr{R}$ and $h^{(2k)}_{\alpha,\kappa_\sigma}(0) = \kappa_\sigma(2\pi)^{2k} + \kappa_\sigma(\kappa_\sigma - 1)(2\alpha)^{2k} - (\kappa_\sigma - 1)(2(\pi - \alpha))^{2k} \leq 0$ for $k$ large enough. Moreover, one can check that if the property $h^{(2k)}_{\alpha,\kappa_\sigma}(0) \leq 0$ is verified for some $k \geq 2$, then there holds $h^{(2m)}_{\alpha,\kappa_\sigma}(0) \leq 0$ for all $m \geq k$. Since we have $h^{(2k+1)}_{\alpha,\kappa_\sigma}(0) = 0$ for all $k \in \mathbb{N}$, we deduce that there hold $h^{(m)}_{\alpha,\kappa_\sigma}(0) \geq 0$ for $m \in \{1, \dots, 2k_0\}$ and $h^{(m)}_{\alpha,\kappa_\sigma}(0) \leq 0$ for $m \in \{2k_0 - 1, 2k_0, 2k_0 + 1, \dots\}$. It is easy to prove that for $m$ large enough, $\eta \mapsto h^{(m)}_{\alpha,\kappa_\sigma}(\eta)$ is strictly negative on $(0; +\infty)$. This allows to show that $\eta \mapsto h^{(m)}_{\alpha,\kappa_\sigma}(\eta)$ is strictly negative on $(0; +\infty)$ for $m \geq 2k_0 - 1$. Since $h^{(2k_0)}_{\alpha,\kappa_\sigma}(0) > 0$, the continuous function $\eta \mapsto h^{(2k_0)}_{\alpha,\kappa_\sigma}(\eta)$ is strictly positive for small $\eta$, strictly negative for large $\eta$, and strictly decreasing. Therefore, it vanishes exactly once on $(0; +\infty)$, for some $\eta_{2k_0}$. The continuous function $\eta \mapsto h^{(2k_0-1)}_{\alpha,\kappa_\sigma}(\eta)$ satisfies $h^{(2k_0-1)}_{\alpha,\kappa_\sigma}(0) = 0$, it is strictly increasing on $(0; \eta_{2k_0})$, strictly decreasing on $(\eta_{2k_0}; +\infty)$ and tends to $-\infty$ when $\eta \to +\infty$. Consequently, it vanishes exactly once on $(0; +\infty)$. Repeating the process we prove by induction that $\eta \mapsto h_{\alpha,\kappa_\sigma}(\eta)$ vanishes exactly once on $(0; +\infty)$. $\square$

**Proposition 5.2.** *Assume that $(\alpha, \kappa_\sigma) \in ((0;\pi) \times (-\infty;0)) \setminus \overline{\mathscr{R}}$ where $\mathscr{R}$ is the region defined in (39). Then, the set of singular exponents $\Lambda_{\alpha,\kappa_\sigma}$ is such that $\Lambda_{\alpha,\kappa_\sigma} \cap \{\lambda \in \mathbb{C} \setminus \{1\} \,|\, \Re e\, \lambda = 1\} = \emptyset$.*

*Proof.* When $(\alpha, \kappa_\sigma) \in ((0;\pi) \times (-\infty;0)) \setminus \overline{\mathscr{R}}$, there holds $h^{(2)}_{\alpha,\kappa_\sigma}(0) = g_\alpha(\kappa_\sigma) < 0$. Our goal is to prove that when $(\alpha, \kappa_\sigma) \in ((0;\pi) \times (-\infty;0)) \setminus \overline{\mathscr{R}}$, we have actually $h^{(m)}_{\alpha,\kappa_\sigma}(0) \leq 0$ for all $m \in \mathbb{N}$. Indeed, since for $m$ large enough, $\eta \mapsto h^{(m)}_{\alpha,\kappa_\sigma}(\eta)$ is strictly negative on $(0;+\infty)$, this will allow to prove by induction that $\eta \mapsto h_{\alpha,\kappa_\sigma}(\eta)$ is strictly negative on $(0;+\infty)$.

First, it is easy to see that there holds $h^{(2k+1)}_{\alpha,\kappa_\sigma}(0) = 0$ for all $k \in \mathbb{N}$. The property $h_{\alpha,\kappa_\sigma}(0) = 0$ is also simple to obtain. Let us show that $h^{(4)}_{\alpha,\kappa_\sigma}(0) \leq 0$. This will end the proof because for $k \geq 2$, we notice that if $h^{(2k)}_{\alpha,\kappa_\sigma}(0) = \kappa_\sigma (2\pi)^{2k} + \kappa_\sigma(\kappa_\sigma - 1)(2\alpha)^{2k} - (\kappa_\sigma - 1)(2(\pi - \alpha))^{2k} \leq 0$, then $h^{(2(k+1))}_{\alpha,\kappa_\sigma}(0) \leq 0$.

We compute $h^{(4)}_{\alpha,\kappa_\sigma}(0) = \tilde{g}_\alpha(\kappa_\sigma)$ with, for $\kappa_\sigma \in (-\infty;0)$,

$$\tilde{g}_\alpha(\kappa_\sigma) = \kappa_\sigma(2\pi)^4 + \kappa_\sigma(\kappa_\sigma - 1)(2\alpha)^4 - (\kappa_\sigma - 1)(2(\pi - \alpha))^4.$$

For a given angle $\alpha \in (0;\pi)$, the function $\kappa_\sigma \mapsto \tilde{g}_\alpha(\kappa_\sigma)$ is a polynomial of degree two. By a simple study, we find that for all $\alpha \in (0;\pi)$, the equation $\tilde{g}_\alpha(\kappa_\sigma) = 0$ has exactly two roots $\tilde{\ell}_-(\alpha) < \tilde{\ell}_+(\alpha) < 0$ on $(-\infty;0)$. Moreover, $\tilde{g}_\alpha$ is negative on $(\tilde{\ell}_-(\alpha); \tilde{\ell}_+(\alpha))$. Therefore, to establish that $h^{(4)}_{\alpha,\kappa_\sigma}(0) = \tilde{g}_\alpha(\kappa_\sigma) \leq 0$ when $(\alpha, \kappa_\sigma) \in ((0;\pi) \times (-\infty;0)) \setminus \overline{\mathscr{R}}$, it remains to prove that there holds

$$\tilde{\ell}_-(\alpha) \leq \ell_-(\alpha) = -\frac{\pi - \alpha + \sin(\pi - \alpha)}{\alpha - \sin \alpha} \leq \ell_+(\alpha) = -\frac{\pi - \alpha - \sin(\pi - \alpha)}{\alpha + \sin \alpha} \leq \tilde{\ell}_+(\alpha).$$

To obtain this result, we show that $\tilde{g}_\alpha(\ell_-(\alpha)) \leq 0$ and $\tilde{g}_\alpha(\ell_+(\alpha)) \leq 0$. We find

$$\frac{\alpha - \sin \alpha}{\pi - \alpha + \sin(\pi - \alpha)} \tilde{g}_\alpha(\ell_-(\alpha)) = (2\pi)^4 \gamma(\alpha)$$

with

$$\gamma(\alpha) = -1 + \frac{\pi}{\alpha - \sin \alpha} \frac{\alpha^4}{\pi^4} + \frac{\pi}{\pi - \alpha + \sin(\pi - \alpha)} \frac{(\pi - \alpha)^4}{\pi^4}.$$

Using the inequalities $\alpha - \sin \alpha \geq \alpha^3/\pi^2$ and $\sin(\pi - \alpha) \geq 0$ for $\alpha \in (0;\pi)$, we obtain

$$\gamma(\alpha) \leq -1 + \frac{\alpha}{\pi} + \frac{(\pi - \alpha)^3}{\pi^3} = -\alpha(\pi - \alpha)(2\pi - \alpha)/\pi^3. \tag{51}$$

Clearly, the right hand side of (51) is negative on $(0;\pi)$. Thus, there holds $\tilde{g}_\alpha(\ell_-(\alpha)) \leq 0$. Noticing that

$$\frac{\alpha + \sin \alpha}{\pi - \alpha - \sin(\pi - \alpha)} \tilde{g}_\alpha(\ell_+(\alpha)) = (2\pi)^4 \gamma(\pi - \alpha),$$

we deduce that we also have $\tilde{g}_\alpha(\ell_+(\alpha)) \leq 0$. This concludes the proof. $\square$

## Acknowledgments

Part of this work has been realized while the author was visiting the Department of Mathematical Sciences of the University of Delaware thanks to the NSF Grant DMS-1106972. The support of the Associate Team ISIP between DeFi-INRIA at École Polytechnique and the University of Delaware is gratefully acknowledged. The author also wishes to express its gratitude to Fioralba Cakoni and David Colton for their very warm welcome, and to Jérémy Firozaly for his help concerning the computation of singularities of §3.2.1 during his master thesis.